\documentclass{article}
\usepackage{amssymb}

\usepackage{amsmath}

\input{tcilatex}

\begin{document}

\title{\textbf{A CENTRAL LIMIT THEOREM AND HIGHER ORDER RESULTS FOR THE
ANGULAR BISPECTRUM\thanks{%
I am very much grateful to M.W.Baldoni for many discussions and explanations
on the role of Clebsch-Gordan coefficients in representation theory. Usual
disclaimers apply.}}}
\author{Domenico Marinucci \\
Department of Mathematics, University of Rome ``Tor Vergata''\\
marinucc@mat.uniroma2.it}
\maketitle

\begin{abstract}
The angular bispectrum of spherical random fields has recently gained an
enormous importance, especially in connection with statistical inference on
cosmological data. In this paper, we provide expressions for its moments of
arbitrary order and we use these results to establish a multivariate central
limit theorem and higher order approximations. The results rely upon
combinatorial methods from graph theory and a detailed investigation for the
asymptotic behaviour of Clebsch-Gordan coefficients; the latter are widely
used in representation theory and quantum theory of angular momentum.

\begin{itemize}
\item \textbf{AMS\ 2000 Classification:} Primary 60G60; Secondary 60F05,
62M15, 62M40

\item \textbf{Key words and Phrases: }Spherical Random Fields, Angular
Bispectrum, Graphs, Clebsch-Gordan Coeffficients, Central Limit Theorem,
Higher Order Approximations.
\end{itemize}
\end{abstract}

\begin{center}
\textbf{1. INTRODUCTION}

\qquad
\end{center}

Let $T(\theta ,\varphi )$ be a random field indexed by the unit sphere $%
S^{2} $, i.e. $0\leq \theta \leq \pi $ and $0\leq \varphi <2\pi $. We assume
that $T(\theta ,\varphi )$ has zero mean, finite variance and it is mean
square continuous and isotropic, i.e. its covariance is invariant with
respect to the group of rotations. For isotropic fields, the following
spectral representation holds in mean square sense (Yaglom (1986), Leonenko
(1999)): 
\begin{equation}
T(\theta ,\varphi )=\sum_{l=1}^{\infty }\sum_{m=-l}^{l}a_{lm}Y_{lm}(\theta
,\varphi )\text{ . }  \label{specrap}
\end{equation}%
Here, we have introduced the spherical harmonics (see Varshalovich, Moskalev
and Khersonskii (hereafter VMK) (1988), chapter 5), defined by 
\begin{align*}
Y_{lm}(\theta ,\varphi )& :=\sqrt{\frac{2l+1}{4\pi }\frac{(l-m)!}{(l+m)!}}%
P_{lm}(\cos \theta )\exp (im\varphi )\text{ , for }m\geq 0\text{ ,} \\
Y_{lm}(\theta ,\varphi )& :=(-1)^{m}Y_{l,-m}^{\ast }(\theta ,\varphi )\text{
, for }m<0\text{ },
\end{align*}%
where the asterisk denotes complex conjugation and $P_{lm}(\cos \theta )$
denotes the associated Legendre polynomial of degree $l,m,$ i.e. 
\begin{align*}
P_{lm}(x)& :=(-1)^{m}(1-x^{2})^{m/2}\frac{d^{m}}{dx^{m}}P_{l}(x)\text{ , }%
P_{l}(x)=\frac{1}{2^{l}l!}\frac{d^{l}}{dx^{l}}(x^{2}-1)^{l}, \\
m& =0,1,2,...,l\text{ , }l=1,2,3,....\text{ .}
\end{align*}%
The triangular array $\left\{ a_{lm}\right\} $ represents a set of random
coefficients, which can be obtained from $T(\theta ,\varphi )$ through the
inversion formula 
\begin{equation}
a_{lm}=\int_{-\pi }^{\pi }\int_{0}^{\pi }T(\theta ,\varphi )Y_{lm}^{\ast
}(\theta ,\varphi )\sin \theta d\theta d\varphi \text{ , }m=0,\pm 1,...,\pm l%
\text{ },\text{ }l=1,2,...\text{ ;}  \label{alm}
\end{equation}%
see for instance Kim and Koo (2002), Kim, Koo and Park (2004) for a review
of Fourier analysis on $S^{2}$. The coefficients $a_{lm}$ are
complex-valued, zero-mean and uncorrelated; hence, if $T(\theta ,\varphi )$
is Gaussian, they have a complex Gaussian distribution, and they are
independent over $l$ and $m\geq 0$ (although $a_{l,-m}=(-1)^{m}a_{lm}^{\ast
} $), with variance $E|a_{lm}|^{2}=C_{l},$ $m=0,\pm 1,...,\pm l.$ The index $%
l$ is usually labeled a multipole; approximately, each multipole corresponds
to an angular resolution given by $\pi /l$.

The sequence $\left\{ C_{l}\right\} $ denotes the angular power spectrum: we
shall always assume that $C_{l}$ is strictly positive, for all values of $l.$
As well-known, if the field is Gaussian the angular power spectrum
completely identifies its dependence structure. For non-Gaussian fields, the
dependence structure becomes much richer, and higher order moments of the $%
a_{lm}$'s are of interest; this leads to the analysis of so-called higher
order angular power spectra.

The analysis of spherical random fields has recently gained momentum, due to
strong empirical motivations arising especially (but not exclusively) from
cosmology and astrophysics. In particular, an enormous attention has been
drawn by issues connected with the statistical analysis of Cosmic Microwave
Background radiation (CMB). CMB can be viewed as a snapshot of the Universe
approximately $3\times 10^{5}$ years after the Big Bang (Peebles (1993),
Peacock (1999)). A number of experiments are aimed at measuring this
radiation: we mention in particular two satellite missions, namely WMAP by
NASA, which released the first full-sky of CMB fluctuations in February
2003, with much more detailed data to come in the years to come, and Planck
by ESA, which is due to be launched in Spring 2007 and expected to provide
maps with much greater resolution. Over the next ten years, an immense
amount of cosmological information is expected from these huge data sets; at
the same time, the analysis of such data sets posits a remarkable challenge
to statistical methodology. In particular, several papers have focussed on
testing for non-Gaussianity by a variety of nonparametric methods (to
mention a few, Dor\'{e}, Colombi and Bouchet (2003), Hansen, Marinucci and
Vittorio (2003), Park (2004), Marinucci and Piccioni (2004), Jin et al.
(2005)). The majority of efforts has focussed on the angular bispectrum,
which is considered an optimal statistic to verify the accuracy of the
so-called inflationary scenario, the leading paradigm for the dynamics of
the Big Bang. See for instance Phillips and Kogut (2000), Komatsu and
Spergel (2001), Bartolo, Matarrese and Riotto (2002), Komatsu et al.
(2002,2003), Babich (2005) and several others; a review is in Marinucci
(2004).

The angular bispectrum can be viewed as the harmonic transform of the
three-point angular correlation function, whereas the angular power spectrum
is the Legendre transform of the (two-point) angular correlation function.
Write $\Omega _{i}=(\theta _{i},\varphi _{i}),$ for $i=1,2,3;$ we have 
\begin{equation}
ET(\Omega _{1})T(\Omega _{2})T(\Omega
_{3})=\sum_{l_{1},l_{2},l_{3}=1}^{\infty
}\sum_{m_{1},m_{2},m_{3}}B_{l_{1}l_{2}l_{3}}^{m_{1}m_{2}m_{3}}Y_{l_{1}m_{1}}(\Omega _{1})Y_{l_{2}m_{2}}(\Omega _{2})Y_{l_{3}m_{3}}(\Omega _{3})%
\text{ ,}  \label{vene1}
\end{equation}%
where the bispectrum $B_{l_{1}l_{2}l_{3}}^{m_{1}m_{2}m_{3}}$ is given by%
\begin{equation}
B_{l_{1}l_{2}l_{3}}^{m_{1}m_{2}m_{3}}=E(a_{l_{1}m_{1}}a_{l_{2}m_{2}}a_{l_{3}m_{3}})%
\text{ .}  \label{vene2}
\end{equation}%
Here, and in the sequel, the sums over $m_{i}$ run from $-l_{i}$ to $l_{i}$,
unless otherwise indicated. Both (\ref{vene1}) and (\ref{vene2}) are clearly
equal to zero for zero-mean Gaussian fields. Moreover, the assumption that
the CMB\ random field is statistically isotropic entails that the right- and
left-hand sides of (\ref{vene1}) should be left unaltered by a rotation of
the coordinate system. Therefore $B_{l_{1}l_{2}l_{3}}^{m_{1}m_{2}m_{3}}$
must take values ensuring that the three-point correlation function on the
left-hand side of (\ref{vene1}) remains unchanged if the three directions $%
\Omega _{1},\Omega _{2}\ $\ and $\Omega _{3}$ are rotated by the same angle.
Careful choices of the orientations entail that the angular bispectrum of an
isotropic field can be non-zero only if $l_{i}\leq l_{j}+l_{k}$ for all
choices of $i,j,k=1,2,3;$ $l_{1}+l_{2}+l_{3}$ is even; and $%
m_{1}+m_{2}+m_{3}=0.$ More generally, Hu (2001) shows that a necessary and
sufficient condition for $B_{l_{1}l_{2}l_{3}}^{m_{1}m_{2}m_{3}}$ to
represent the angular bispectrum of an isotropic random field is that there
exist a real symmetric function of $l_{1},l_{2},l_{3},$ which we denote $%
b_{l_{1}l_{2}l_{3}},$ such that we have the identity%
\begin{equation}
B_{l_{1}l_{2}l_{3}}^{m_{1}m_{2}m_{3}}=\mathcal{G}%
_{l_{1}l_{2}l_{3}}^{m_{1}m_{2}m_{3}}b_{l_{1}l_{2}l_{3}}\text{ ;}
\label{veneid}
\end{equation}%
$b_{l_{1}l_{2}l_{3}}$ is labeled the reduced bispectrum. In (\ref{veneid})
we are using the Gaunt integral $\mathcal{G}%
_{l_{1}l_{2}l_{3}}^{m_{1}m_{2}m_{3}},$ defined by%
\begin{equation*}
\mathcal{G}_{l_{1}l_{2}l_{3}}^{m_{1}m_{2}m_{3}}:=\int_{0}^{\pi
}\int_{0}^{2\pi }Y_{l_{1}m_{1}}(\theta ,\varphi )Y_{l_{2}m_{2}}(\theta
,\varphi )Y_{l_{3}m_{3}}(\theta ,\varphi )\sin \theta d\varphi d\theta
\end{equation*}%
\begin{equation*}
=\left( \frac{(2l_{1}+1)(2l_{2}+1)(2l_{3}+1)}{4\pi }\right) ^{1/2}\left( 
\begin{tabular}{lll}
$l_{1}$ & $l_{2}$ & $l_{3}$ \\ 
$0$ & $0$ & $0$%
\end{tabular}%
\right) \left( 
\begin{tabular}{lll}
$l_{1}$ & $l_{2}$ & $l_{3}$ \\ 
$m_{1}$ & $m_{2}$ & $m_{3}$%
\end{tabular}%
\right) \text{ ,}
\end{equation*}%
where the so-called ``Wigner's $3j$ symbols'' appearing on the second line
are defined by (VMK, expression 8.2.1.5) 
\begin{align*}
\left( 
\begin{array}{ccc}
l_{1} & l_{2} & l_{3} \\ 
m_{1} & m_{2} & m_{3}%
\end{array}%
\right) & :=(-1)^{l_{3}+m_{3}+l_{2}+m_{2}}\left[ \frac{%
(l_{1}+l_{2}-l_{3})!(l_{1}-l_{2}+l_{3})!(l_{1}-l_{2}+l_{3})!}{%
(l_{1}+l_{2}+l_{3}+1)!}\right] ^{1/2} \\
& \times \left[ \frac{(l_{3}+m_{3})!(l_{3}-m_{3})!}{%
(l_{1}+m_{1})!(l_{1}-m_{1})!(l_{2}+m_{2})!(l_{2}-m_{2})!}\right] ^{1/2} \\
& \times \sum_{z}\frac{(-1)^{z}(l_{2}+l_{3}+m_{1}-z)!(l_{1}-m_{1}+z)!}{%
z!(l_{2}+l_{3}-l_{1}-z)!(l_{3}+m_{3}-z)!(l_{1}-l_{2}-m_{3}+z)!}\text{ ,}
\end{align*}%
where the summation runs over all $z$'s such that the factorials are
non-negative. Note that the Wigner's $3j$ are invariant with respect to
permutations of the pairs $(l_{i},m_{i}).$

In view of (\ref{veneid}), the dependence of the bispectrum ordinates on $%
m_{1},m_{2},m_{3}$ does not carry any physical information if the field is
isotropic; hence it can be eliminated by focussing on the angular averaged
bispectrum, defined by%
\begin{align}
B_{l_{1}l_{2}l_{3}}&
:=\sum_{m_{1}=-l_{1}}^{l_{1}}\sum_{m_{2}=-l_{2}}^{l_{2}}%
\sum_{m_{3}=-l_{3}}^{l_{3}}\left( 
\begin{tabular}{lll}
$l_{1}$ & $l_{2}$ & $l_{3}$ \\ 
$m_{1}$ & $m_{2}$ & $m_{3}$%
\end{tabular}%
\right) B_{l_{1}l_{2}l_{3}}^{m_{1}m_{2}m_{3}}  \notag \\
& =\left( \frac{(2l_{1}+1)(2l_{2}+1)(2l_{3}+1)}{4\pi }\right) ^{1/2}\left( 
\begin{tabular}{lll}
$l_{1}$ & $l_{2}$ & $l_{3}$ \\ 
$0$ & $0$ & $0$%
\end{tabular}%
\right) b_{l_{1}l_{2}l_{3}}\text{ ,}  \label{parr}
\end{align}%
where we have used the orthogonality condition%
\begin{equation*}
\sum_{m_{1}=-l_{1}}^{l_{1}}\sum_{m_{2}=-l_{2}}^{l_{2}}%
\sum_{m_{3}=-l_{3}}^{l_{3}}\left( 
\begin{tabular}{lll}
$l_{1}$ & $l_{2}$ & $l_{3}$ \\ 
$m_{1}$ & $m_{2}$ & $m_{3}$%
\end{tabular}%
\right) ^{2}=1\text{ .}
\end{equation*}%
The minimum mean square error estimator of the bispectrum is provided by (Hu
2001) 
\begin{equation*}
\widehat{B}_{l_{1}l_{2}l_{3}}:=\sum_{m_{1}=-l_{1}}^{l_{1}}%
\sum_{m_{2}=-l_{2}}^{l_{2}}\sum_{m_{3}=-l_{3}}^{l_{3}}\left( 
\begin{tabular}{lll}
$l_{1}$ & $l_{2}$ & $l_{3}$ \\ 
$m_{1}$ & $m_{2}$ & $m_{3}$%
\end{tabular}%
\right) (a_{l_{1}m_{1}}a_{l_{2}m_{2}}a_{l_{3}m_{3}})\text{ .}
\end{equation*}%
The statistic $\widehat{B}_{l_{1}l_{2}l_{3}}$ is called the (sample) angle
averaged bispectrum; for any realization of the random field $T,$ it is a
real-valued scalar, which does not depend on the choice of the coordinate
axes and it is invariant with respect to permutation of its arguments $%
l_{1},l_{2},l_{3}.$

Under Gaussianity, the bispectrum can be easily made model-independent,
namely we can focus on the normalized bispectrum, which we define by%
\begin{equation}
I_{l_{1}l_{2}l_{3}}:=(-1)^{(l_{1}+l_{2}+l_{3})/2}\frac{\widehat{B}%
_{l_{1}l_{2}l_{3}}}{\sqrt{C_{l_{1}}C_{l_{2}}C_{l_{3}}}}\text{ .}
\label{ult1}
\end{equation}%
The factor $(-1)^{(l_{1}+l_{2}+l_{3})/2}$ is usually not included in the
definition of the normalized bispectrum; it corresponds, however, to the
sign of the Wigner's coefficients for $m_{1}=m_{2}=m_{3}=0$, and thus it
seems natural to include it to ensure that $I_{l_{1}l_{2}l_{3}}$ and $%
b_{l_{1}l_{2}l_{3}}$ share the same parity (see (\ref{parr})).

In practice $I_{l_{1}l_{2}l_{3}}$ is unfeasible because $C_{l}$ is unknown.
A natural estimator for $C_{l}$ is 
\begin{equation}
\widehat{C}_{l}:=\frac{1}{2l+1}\sum_{m=-l}^{l}|a_{lm}|^{2}\text{ , }l=1,2,...%
\text{ ,}  \label{clest}
\end{equation}%
which is clearly unbiased (see also Arjunwadkar et al. (2004)). Thus $%
I_{l_{1}l_{2}l_{3}}$ can be replaced by the feasible statistic 
\begin{equation*}
\widehat{I}_{l_{1}l_{2}l_{3}}:=(-1)^{(l_{1}+l_{2}+l_{3})/2}\frac{\widehat{B}%
_{l_{1}l_{2}l_{3}}}{\sqrt{\widehat{C}_{l_{1}}\widehat{C}_{l_{2}}\widehat{C}%
_{l_{3}}}}\text{ .}
\end{equation*}

Although the bispectrum has been the object of an enormous attention in the
cosmological literature, very few analytic results are so far available on
its probabilistic properties. In a previous paper (Marinucci (2005)), we
established some bounds on the behaviour of its first eight moments, and we
used these results to establish the asymptotic behaviour of some functionals
of the bispectrum array; such functionals were proposed to build
nonparametric tests for non-Gaussianity. The results on these moments where
established by means of a direct analysis of cross-products of Wigner's $3j$
coefficients; this analysis was performed by means of explicit summation
formulae originated from the quantum angular momentum literature (see VMK\
for a very detailed collection of results). In the present paper, these
results are made sharper and extended to moments of arbitrary orders by
means of a more general argument. More precisely, we show how it is possible
to associate to higher moments the coefficients of unitary matrices
transforming alternative bases of tensor product spaces generated by the
spherical harmonics. We are then able to obtain combinatorial expressions
for (cross-)moments of arbitrary orders. These results are then exploited to
obtain multivariate central limit theorems and higher order approximations.
It should be noted that the asymptotic theory presented in this work is of a
fixed-domain type, a framework which has become increasingly popular in
recent years, see for instance Stein (1999) or Loh (2005).

The structure of the paper is as follows: in Section 2 we review some basic
combinatorial material on diagrams and graphs; in Section 3 we present our
general results on moments and we exploit it to obtain multivariate central
limit theorems for the angular bispectrum with known or unknown $C_{l};$ in
Section 4, we discuss higher-order approximations.

\qquad

\begin{center}
\textbf{2. DIAGRAMS AND GRAPHS}
\end{center}

\qquad

We shall review here some elementary notions from graph theory, which is
widely used in physics when handling Wigner's $3j$ coefficients (see VMK,
chapter 11). Take $i=1,...,p$ and $j=1,...,q_{i},$ and consider the set of
indexes:%
\begin{equation*}
T=\left\{ 
\begin{array}{cccc}
(1,1) & ... & .... & (1,q_{1}) \\ 
... & ... & ... & ... \\ 
(p,1) & ... & ... & (p,q_{p})%
\end{array}%
\right\} \text{ ;}
\end{equation*}%
we stress that the number of columns $q_{i}$ need not be the same for each
row $i.$ A \emph{diagram} $\gamma $ is any partition of the elements of $T$
into pairs like $\left\{ (i_{1},j_{1}),(i_{2},j_{2})\right\} :$ these pairs
are called the \emph{edges} of the diagram. We label $\Gamma (T)$ the family
of these diagrams$.$ We also note that if we identify each row $i_{k}$ with
a \emph{vertex} (or \emph{node)}, and view these vertexes as linked together
by the edges $\left\{ (i_{k},j_{k}),(i_{k^{\prime }},j_{k^{\prime
}})\right\} =i_{k}i_{k^{\prime }},$ then it is possible to associate to each
diagram a \emph{graph. }As it is well known, a graph is an ordered pair $%
(I,E)$ where $I$ is non-empty (in our case the set of the rows of the
diagram), and $E$ is a set of unordered pairs of vertexes (in our cases, the
pair of rows that are linked in a diagram). We consider only graphs which
are not directed, that is, $(i_{1}i_{2})$ and $(i_{2}i_{1})$ identify the
same edge; however, we do allow for repetitions of edges (two rows may be
linked twice), in which case the term \emph{multigraph}$\mathbb{\ }$is more
appropriate. A graph carries less information than a diagram (the
information on the ``columns'', i.e. the second element $j_{k},$ is
neglected) but it is much easier to represent pictorially. We shall use some
result on graphs below; with a slight abuse of notation, we denote the graph 
$\gamma $ with the same letter as the corresponding diagram.

We say that

a) A diagram has a \emph{flat edge} if there is at least a pair $\left\{
(i_{1},j_{1}),(i_{2},j_{2})\right\} $ such that $i_{1}=i_{2};$ we write $%
\gamma \in \Gamma _{F}(T)$ for a diagram with at least a flat edge, and $%
\gamma \in \Gamma _{\overline{F}}(T)$ otherwise. A graph corresponding to a
diagram with a flat edge includes an edge $i_{k}i_{k}$ which arrives in the
same vertex where it started; for these circumstances the term \emph{%
pseudograph} is preferred by some authors (e.g. Foulds (1992)).

b) A diagram $\gamma \in \Gamma _{\overline{F}}(T)$ is \emph{connected} if
it is not possible to partition the $i$'s into two sets $A,B$ such that
there are no edges with $i_{1}\in A$ and $i_{2}\in B.$ We write $\gamma \in
\Gamma _{C}(T)$ for connected diagrams, $\gamma \in \Gamma _{\overline{C}%
}(T) $ otherwise. Obviously a diagram is connected if and only if the
corresponding graph is connected, in the standard sense.

c) A diagram $\gamma \in \Gamma _{\overline{F}}(T)$ is \emph{paired} if,
considering any two set of edges $\left\{
(i_{1},j_{1}),(i_{2},j_{2})\right\} $ and $\left\{
(i_{3},j_{3}),(i_{4},j_{4})\right\} ,$ then $i_{1}=i_{3}$ implies $%
i_{2}=i_{4};$ in words, the rows are completely coupled two by two. We write 
$\gamma \in \Gamma _{P}(T)$ for paired diagrams.\FRAME{dtbpFUX}{3in}{2.0003in%
}{0pt}{\Qcb{Figure I: $\protect\gamma \in \Gamma _{P}(T)$}}{}{Plot}{\special%
{language "Scientific Word";type "MAPLEPLOT";width 3in;height 2.0003in;depth
0pt;display "USEDEF";plot_snapshots TRUE;mustRecompute FALSE;lastEngine
"Maple";xmin "1";xmax "2";xviewmin "0.98";xviewmax "2.0204";yviewmin
"0.830000190031882";yviewmax "3.17089980982194";plottype 4;constrained
TRUE;numpoints 49;plotstyle "patch";axesstyle "none";xis \TEXUX{x};var1name
\TEXUX{$x$};function
\TEXUX{\EQN{6}{1}{}{}{\RD{%
\CELL{((x-1)^{1.5}(2-x)^{1.5}+1)}}{1}{}{}{}}};linecolor "black";linestyle
1;pointstyle "point";linethickness 1;lineAttributes "Solid";var1range
"1,2";num-x-gridlines 49;curveColor "[flat::RGB:0000000000]";curveStyle
"Line";rangeset"X";function \TEXUX{$(-(x-1)^{1.5}(2-x)^{1.5}+1)$};linecolor
"black";linestyle 1;pointstyle "point";linethickness 1;lineAttributes
"Solid";var1range "1,2";num-x-gridlines 49;curveColor
"[flat::RGB:0000000000]";curveStyle "Line";rangeset"X";function
\TEXUX{$((x-1)^{1.5}(2-x)^{1.5}+2)$};linecolor "black";linestyle
1;pointstyle "point";linethickness 1;lineAttributes "Solid";var1range
"1,2";num-x-gridlines 49;curveColor "[flat::RGB:0000000000]";curveStyle
"Line";rangeset"X";function \TEXUX{$(-(x-1)^{1.5}(2-x)^{1.5}+2)$};linecolor
"black";linestyle 1;pointstyle "point";linethickness 1;lineAttributes
"Solid";var1range "1,2";num-x-gridlines 49;curveColor
"[flat::RGB:0000000000]";curveStyle "Line";rangeset"X";function
\TEXUX{$\left[ 1,1,2,1\right] $};linecolor "black";linestyle 1;pointstyle
"point";linethickness 1;lineAttributes "Solid";curveColor
"[flat::RGB:0000000000]";curveStyle "Line";function \TEXUX{$\left[
1,2,2,2\right] $};linecolor "black";linestyle 1;pointstyle
"point";linethickness 1;lineAttributes "Solid";curveColor
"[flat::RGB:0000000000]";curveStyle "Line";function
\TEXUX{$(-(x-1)^{1.5}(2-x)^{1.5}+3)$};linecolor "black";linestyle
1;pointstyle "point";linethickness 1;lineAttributes "Solid";var1range
"1,2";num-x-gridlines 49;curveColor "[flat::RGB:0000000000]";curveStyle
"Line";rangeset"X";function \TEXUX{$((x-1)^{1.5}(2-x)^{1.5}+2)$};linecolor
"black";linestyle 1;pointstyle "point";linethickness 1;lineAttributes
"Solid";var1range "1,2";num-x-gridlines 49;curveColor
"[flat::RGB:0000000000]";curveStyle "Line";rangeset"X";function
\TEXUX{$\left[ 1,3,2,3\right] $};linecolor "black";linestyle 1;pointstyle
"point";linethickness 1;lineAttributes "Solid";curveColor
"[flat::RGB:0000000000]";curveStyle "Line";rangeset"X";function
\TEXUX{$((x-1)^{1.5}(2-x)^{1.5}+3)$};linecolor "black";linestyle
1;pointstyle "point";linethickness 1;lineAttributes "Solid";var1range
"1,2";num-x-gridlines 49;curveColor "[flat::RGB:0000000000]";curveStyle
"Line";rangeset"X";function \TEXUX{$\left[ 1,1,2,1\right] $};linecolor
"black";linestyle 1;pointplot TRUE;pointstyle "box";linethickness
1;lineAttributes "Solid";curveColor "[flat::RGB:0000000000]";curveStyle
"Point";function \TEXUX{$\left[ 1,2,2,2\right] $};linecolor
"black";linestyle 1;pointplot TRUE;pointstyle "box";linethickness
1;lineAttributes "Solid";curveColor "[flat::RGB:0000000000]";curveStyle
"Point";function \TEXUX{$\left[ 1,3,2,3\right] $};linecolor
"black";linestyle 1;pointplot TRUE;pointstyle "box";linethickness
1;lineAttributes "Solid";curveColor "[flat::RGB:0000000000]";curveStyle
"Point";rangeset"X";valid_file "T";tempfilename
'IIDRFK01.wmf';tempfile-properties "XPR";}}

d) We shall say a diagram has a \emph{k-loop} if there exist a sequence of $%
k $ edges 
\begin{equation*}
\left\{ (i_{1},j_{1}),(i_{2},j_{2})\right\} ,...,\left\{
(i_{k},j_{k}),(i_{k+1},j_{k+1})\right\} =(i_{1}i_{2}),...,(i_{k}i_{k+1})
\end{equation*}%
such that $i_{1}=i_{k+1};$ we write $\gamma \in \Gamma _{L(k)}(T)$ for
diagrams with a $k$-loop and no loop of order smaller than $k.$

Note that $\Gamma _{F}(T)=\Gamma _{L(1)}(T)$ (a flat edge is a 1-loop);
also, we write 
\begin{equation*}
\Gamma _{CL(k)}(T)=\Gamma _{C}(T)\cap \Gamma _{L(k)}(T)\text{ }
\end{equation*}%
for connected diagrams with k-loops, and $\Gamma _{C\overline{L(k)}}(T)$ for
connected diagrams with no loops of order $k$ or smaller. For instance, a
connected diagram belongs to $\Gamma _{C\overline{L(2)}}(T)$ if there are
neither flat edges nor two edges $\left\{
(i_{1},j_{1}),(i_{2},j_{2})\right\} $ and $\left\{
(i_{3},j_{3}),(i_{4},j_{4})\right\} $ such that $i_{1}=i_{3}$ and $%
i_{2}=i_{4};$ in words, there are no pairs of rows which are connected twice.

e) A \emph{tree} is a graph with no loops (written $\gamma \in \Gamma
_{T}(T))$.

Graphs and diagrams play a key role to evaluate the behaviour of moments of
the bispectrum; to this issue we devote the next section.

\ \newline

\begin{center}
\textbf{3.} \textbf{A CENTRAL LIMIT THEOREM FOR THE BISPECTRUM}

\ 
\end{center}

In this section, we shall investigate the behaviour of the higher order
moments for the normalized bispectrum (\ref{ult1}), under the assumption of
Gaussianity; to this aim, we define%
\begin{equation*}
\Delta _{l_{1}l_{2}l_{3}}:=1+\delta _{l_{1}}^{l_{2}}+\delta
_{l_{2}}^{l_{3}}+3\delta _{l_{1}}^{l_{3}}=\left\{ 
\begin{array}{c}
1\text{ for }l_{1}<l_{2}<l_{3} \\ 
2\text{ for }l_{1}=l_{2}<l_{3}\text{ or }l_{1}<l_{2}=l_{3} \\ 
6\text{ for }l_{1}=l_{2}=l_{3}%
\end{array}%
\right. ;
\end{equation*}%
here and in the sequel, $\delta _{a}^{b}$ denotes Kronecker's delta, that is 
$\delta _{a}^{b}=1$ for $a=b,$ zero otherwise.

Under Gaussianity, it is obvious that the expectation of all odd powers of $%
I_{l_{1}l_{2}l_{3}}$ is zero. To analyze the behaviour of even powers, we
first recall that, for a multivariate Gaussian vector $(x_{1},...,x_{2k}),$
we have the following diagram formula%
\begin{equation}
E(x_{1}\times x_{2}\times ...\times x_{2k})=\sum (Ex_{i_{1}}x_{i_{2}})\times
...\times (Ex_{i_{2k-1}}x_{i_{2k}})\text{ ,}  \label{adle}
\end{equation}%
where the sum is over all the $(2k)!/(k!2^{k})$ different ways of grouping $%
(x_{1},...,x_{2k})$ into pairs (see for instance Giraitis and Surgailis
(1987)). Even powers of $I_{l_{1}l_{2}l_{3}}$ yield even powers of the $%
a_{lm}$'s, which have a complex Gaussian distribution, weighted by Wigner's $%
3j$ coefficients.

In the sequel, unless otherwise specified, we rearrange terms so that $%
l_{i1}\leq l_{i2}\leq l_{i3}$ for all $i.$

\qquad \newline
\ \newline
\textbf{Theorem 3.1 }Assume that $(l_{i1},l_{i2},l_{i3})\neq (l_{i^{\prime
}1},l_{i^{\prime }2},l_{i^{\prime }3})$ whenever $i\neq i^{\prime }.$ There
exist an absolute constant $K_{p_{1}...p_{I}}$ such that, for $p_{i}\geq 1$
, $i=1,...,I$%
\begin{equation}
\left| E\left\{ \prod_{i=1}^{I}I_{l_{i1}l_{i2}l_{i3}}^{2p_{i}}\right\}
-\prod_{i=1}^{I}\left\{ (2p_{i}-1)!!\Delta
_{l_{i1}l_{i2}l_{i3}}^{p_{i}}\right\} \right| \leq \frac{K_{p_{1}...p_{I}}}{%
2l_{11}+1}  \label{teo3.1}
\end{equation}%
always, where $(2p-1)!!=(2p-1)\times (2p-3)\times ...\times 1$.

\ \newline
\textbf{Remark 3.1} The condition that $(l_{i1},l_{i2},l_{i3})\neq
(l_{i^{\prime }1},l_{i^{\prime }2},l_{i^{\prime }3})$ whenever $i\neq
i^{\prime }$ is merely notational and entails no loss of generality; indeed,
whenever $(l_{i1},l_{i2},l_{i3})=(l_{i^{\prime }1},l_{i^{\prime
}2},l_{i^{\prime }3})$ it suffices to identify the two indexes and change
the values of $p_{i}$ accordingly.

\textbf{\ }\newline
\textbf{Proof }The proof is lengthy and computationally burdensome, so
before we proceed we find it useful to sketch heuristically its main
features. The first step is to notice that, in view of (\ref{adle}), higher
order moments can be associated with sums over all possible graphs
configurations of cross-products of Wigner's $3j$ coefficients. Our aim
below will be to show that the contribution of each of these components is
determined by its degree of connectivity. Indeed, the leading term will be
provided by paired graphs $\gamma \in \Gamma _{P}(T),$ where the nodes are
partitioned into disjoint pairs. Next to that, we shall show that the
components where at least $p$ nodes are connected are bounded by $O(l_{11}^{-%
\frac{p}{4}})$ (a bound that can be improved for some values of $p).$ This
bound can be obtained by partitioning these connected graphs into trees,
that is subgraphs with no loops, and then associating these trees to the
coefficients of some unitary matrices arising in tensor spaces generated by
spherical harmonics. The proof of (\ref{teo3.1}) can then be simply
concluded by a direct graph-counting argument. Let us now make this argument
rigorous.

We start by introducing some notation, which is to some extent the same as
in Marinucci (2005). We need first to introduce a new set of triples $%
\mathcal{L}=\left\{ \left( \ell _{11},\ell _{12},\ell _{13}\right)
,...\left( \ell _{R1},\ell _{R2},\ell _{R3}\right) \right\} ,$ defined by 
\begin{eqnarray*}
\left( \ell _{r1},\ell _{r2},\ell _{r3}\right) &:&=(l_{11},l_{12},l_{13})%
\text{ for }r=1,...,2p_{1} \\
\left( \ell _{r1},\ell _{r2},\ell _{r3}\right) &:&=(l_{21},l_{22},l_{23})%
\text{ for }r=p_{1}+1,...,2p_{2} \\
&&... \\
\left( \ell _{r1},\ell _{r2},\ell _{r3}\right) &:&=(l_{I1},l_{I2},l_{I3})%
\text{ for }r=p_{I-1}+1,...,2p_{I}\text{ ;}
\end{eqnarray*}%
more explicitly, the set $\mathcal{L}$ is obtained by replicating $2p_{i}$
times each of the $(l_{i1},l_{i2},l_{i3})$ triples. Let $T$ be a set of
indexes $\left\{ (r,k)\right\} $, where $k=1,2,3$ and $r=1,2,...,%
\sum_{i=1}^{I}p_{i};$.for any $\gamma \in \Gamma (T),$ we can define%
\begin{equation}
\delta (\gamma ;\mathcal{L}):=\prod_{\left\{ (r_{u}k_{u}),(r_{u}^{\prime
}k_{u}^{\prime })\right\} \in \gamma }(-1)^{m_{r_{u}k_{u}}}\delta
_{m_{r_{u}k_{u}}}^{-m_{r_{u}^{\prime }k_{u}^{\prime }}}\delta _{\ell
_{r_{u}k_{u}}}^{\ell _{r_{u}k_{u}^{\prime }}}\text{ ;}  \label{invbar2}
\end{equation}%
for brevity's sake, we write $\delta (\gamma )$ rather than $\delta (\gamma ;%
\mathcal{L})$ whenever this causes no ambiguity. Recall that 
\begin{equation}
Ea_{\ell _{r_{1}k_{1}}m_{r_{1}k_{1}}}a_{\ell
_{r_{2}k_{2}}m_{r_{2}k_{2}}}=(-1)^{m_{r_{1}k_{1}}}C_{\ell _{k_{1}}}\delta
_{\ell _{r_{1}k_{1}}}^{\ell _{r_{2}k_{2}}}\delta
_{m_{r_{1}k_{1}}}^{-m_{r_{2}k_{2}}}\text{ }.  \label{invbar3}
\end{equation}%
In view of (\ref{adle}), and because the spherical harmonic coefficients are
(complex) Gaussian distributed, the following formula holds:%
\begin{equation}
E\left\{ \prod_{(r,k)\in T}\frac{a_{\ell _{rk}m_{rk}}}{\sqrt{C_{\ell _{r}}}}%
\right\} =\sum_{\gamma \in \Gamma (T)}\delta (\gamma )\text{ .}  \label{pol1}
\end{equation}%
Write $\left\{ (r,k),.\right\} \in \gamma $ to signify that the pair $%
\left\{ (r,k),(r^{\prime },k^{\prime })\right\} $ belongs to $\gamma ,$ for
some $(r^{\prime },k^{\prime });$ for any diagram $\gamma ,$ we can hence
define%
\begin{equation}
D(\gamma ):=\sum_{\substack{ m_{rk}=-\ell _{rk}  \\ \left\{ (r,k),.\right\}
\in \gamma }}^{\ell _{rk}}\prod_{r:(r,k)\in T}\left( 
\begin{tabular}{lll}
$\ell _{r1}$ & $\ell _{r2}$ & $\ell _{r3}$ \\ 
$m_{r1}$ & $m_{r2}$ & $m_{r3}$%
\end{tabular}%
\right) \delta (\gamma )\text{ .}  \label{invbar7}
\end{equation}%
We define also%
\begin{equation*}
D(A):=\sum_{\gamma \in A}D(\gamma )=\sum_{\gamma \in A}\sum_{\substack{ %
m_{rk}=-\ell _{rk}  \\ \left\{ (r,k),.\right\} \in \gamma }}^{\ell
_{rk}}\prod_{r:(r,k)\in T}\left( 
\begin{tabular}{lll}
$\ell _{1}$ & $\ell _{2}$ & $\ell _{3}$ \\ 
$m_{r1}$ & $m_{r2}$ & $m_{r3}$%
\end{tabular}%
\right) \delta (\gamma )\text{ ;}
\end{equation*}%
in words, $D(.)$ represents the component of the expected value that
corresponds to a particular set of diagrams. Notice that 
\begin{align*}
E\left\{ \prod_{i=1}^{I}I_{l_{1}l_{2}l_{3}}^{2p_{i}}\right\} & =\sum 
_{\substack{ m_{rk}=-\ell _{rk}  \\ r:(r,k)\in T}}^{\ell _{rk}}E\left\{
\prod_{r:(r,k)\in T}\left[ \left( 
\begin{tabular}{lll}
$\ell _{r1}$ & $\ell _{r2}$ & $\ell _{r3}$ \\ 
$m_{r1}$ & $m_{r2}$ & $m_{r3}$%
\end{tabular}%
\right) \prod_{k=1}^{3}\frac{a_{\ell _{rk}m_{rk}}}{\sqrt{C_{\ell _{k}}}}%
\right] \right\} \\
& =\sum_{\gamma \in \Gamma (T)}\sum_{\substack{ m_{rk}=-\ell _{rk}  \\ %
\left\{ (r,k),.\right\} \in \gamma }}^{\ell _{rk}}\prod_{r:(r,k)\in T}\left( 
\begin{tabular}{lll}
$\ell _{r1}$ & $\ell _{r2}$ & $\ell _{r3}$ \\ 
$m_{r1}$ & $m_{r2}$ & $m_{r3}$%
\end{tabular}%
\right) \delta (\gamma ;\mathcal{L}) \\
& =D[\Gamma (T);\mathcal{L}]\text{ .}
\end{align*}%
Now 
\begin{equation*}
D[\Gamma (T);\mathcal{L}]=D[\Gamma _{p}(T);\mathcal{L}]+D[\Gamma
(T)\backslash \Gamma _{P}(T);\mathcal{L}]\text{ ;}
\end{equation*}%
by an identical combinatorial argument as in Marinucci (2005), it is simple
to show that%
\begin{equation*}
D[\Gamma _{p}(T);\mathcal{L}]=\prod_{i=1}^{I}\left\{ (2p_{i}-1)!!\Delta
_{l_{1}l_{2}l_{3}}^{p_{i}}\right\} \text{ .}
\end{equation*}%
To complete the proof, it is then sufficient to establish that 
\begin{equation*}
D[\Gamma (T)\backslash \Gamma _{P}(T);\mathcal{L}]=O\left(
(2l_{11}+1)^{-1}\right) \text{ }.
\end{equation*}%
It is shown in Lemmas 3.1-3.3 in Marinucci (2005) that diagrams with a
1-loop ($\gamma \in \Gamma _{L(1)}(T))$ correspond to summands identically
equal to zero, whereas diagrams with $p$ nodes and loops of orders 2 or 3
can be reduced to terms corresponding to diagrams with $p-2$ nodes times a
factor $O((2l_{11}+1)^{-1}).$\ In the sequel, it is hence sufficient to
focus only on graphs which have no loops of orders 1,2 or 3.

Now call $R$ the set of nodes of the graphs, and partition it into subsets
such that 
\begin{equation*}
R=R_{1}\cup R_{2}\cup ...\cup R_{g}\text{ .}
\end{equation*}%
Then it will also possible to partition $\gamma $ into subdiagrams $\gamma
_{1},\gamma _{2},...,\gamma _{g},$ $\gamma _{12},...,\gamma _{g-1,g}$ such
that $\gamma _{1}$ includes the pairs with both row indexes in $R_{1},$ $%
\gamma _{2}$ includes the pairs with both rows in $R_{2},$ $\gamma _{12}$
includes the pairs with one row in $R_{1}$ and the other in $R_{2},$ and so
on; we assume all internal subdiagrams $\gamma _{i}$ to be non-empty,
whereas this need not be the case for $\gamma _{ij}.$ In terms of edges, $%
\gamma _{1}$ includes the edges that are internal to $R_{1},$ $\gamma _{2}$
includes the edges that are internal to $R_{2},$ $\gamma _{12}$ includes the
edges that connect $R_{1}$ to $R_{2},$ and so forth$.$ Note that, 
\begin{equation*}
\gamma =(\gamma _{1}\cup \gamma _{2}\cup ...\cup \gamma _{g}\cup \gamma
_{12}\cup ...\cup \gamma _{g-1,g})
\end{equation*}%
and we can write%
\begin{eqnarray*}
D(\gamma ) &=&\sum_{i=1}^{g}\sum_{j=i+1}^{g}\sum_{\substack{ m_{rk}=-\ell
_{rk}  \\ \left\{ (r,k),.\right\} \in \gamma _{ij}}}^{\ell _{rk}}\left\{
\prod_{i=1}^{g}\left[ \sum_{\substack{ m_{rk}=-\ell _{rk}  \\ \left\{
(r,k),.\right\} \in \gamma _{i}}}^{\ell _{rk}}\prod_{r\in R_{i}}\left( 
\begin{tabular}{lll}
$\ell _{r1}$ & $\ell _{r2}$ & $\ell _{r3}$ \\ 
$m_{r1}$ & $m_{r2}$ & $m_{r3}$%
\end{tabular}%
\right) \delta (\gamma _{i})\right] \right\} \delta (\gamma _{ij}) \\
&=&\sum_{i=1}^{g}\sum_{j=i+1}^{g}\sum_{\substack{ m_{rk}=-\ell _{rk}  \\ %
\left\{ (r,k),.\right\} \in \gamma _{ij}}}^{\ell _{rk}}\left\{
\prod_{i=1}^{g}X_{R_{i};\gamma _{i}}\right\} \delta (\gamma _{ij})\text{ ,}
\end{eqnarray*}%
where%
\begin{equation*}
X_{R_{i};\gamma _{i}}:=\sum_{\substack{ m_{rk}=-\ell _{rk}  \\ \left\{
(r,k),.\right\} \in \gamma _{i}}}^{\ell _{rk}}\prod_{r\in R_{i}}\left( 
\begin{tabular}{lll}
$\ell _{r1}$ & $\ell _{r2}$ & $\ell _{r3}$ \\ 
$m_{r1}$ & $m_{r2}$ & $m_{r3}$%
\end{tabular}%
\right) \delta (\gamma _{i})\text{ .}
\end{equation*}%
$X_{R_{i};\gamma _{i}}$ can be viewed as a vector whose elements are indexed
by $m_{r_{i},k_{i}}$, where $r_{i}\in R_{i}$ and $\left\{
(r_{i},k_{i}),.\right\} \notin \gamma _{i}$ (indeed those indexes $%
m_{r_{i}k_{i}}$ such that $\left\{ (r_{i},k_{i}),.\right\} \in \gamma _{i}$
have been summed up internally). For instance, for $g=2$ we have%
\begin{eqnarray*}
D(\gamma ) &=&\sum_{\substack{ m_{rk}=-\ell _{rk}  \\ \left\{
(r,k),.\right\} \in \gamma _{12}}}^{\ell _{rk}}\left\{ \sum_{\substack{ %
m_{rk}=-\ell _{rk}  \\ \left\{ (r,k),.\right\} \in \gamma _{1}}}^{\ell
_{rk}}\prod_{r\in R_{1}}\left( 
\begin{tabular}{lll}
$\ell _{r1}$ & $\ell _{r2}$ & $\ell _{r3}$ \\ 
$m_{r1}$ & $m_{r2}$ & $m_{r3}$%
\end{tabular}%
\right) \delta (\gamma _{1})\right\} \\
&&\times \left\{ \sum_{\substack{ m_{rk}=-\ell _{rk}  \\ \left\{
(r,k),.\right\} \in \gamma _{2}}}^{\ell _{rk}}\prod_{r\in R_{2}}\left( 
\begin{tabular}{lll}
$\ell _{r1}$ & $\ell _{r2}$ & $\ell _{r3}$ \\ 
$m_{r1}$ & $m_{r2}$ & $m_{r3}$%
\end{tabular}%
\right) \delta (\gamma _{2})\right\} \delta (\gamma _{12})\text{ .}
\end{eqnarray*}%
In Figure II, we provide a graph with eight nodes $\#(R)=8$ (right), and
then (left) we partition it with $g=2,$ $\#(R_{1})=\#(R_{2})=4;$ the nodes
in $R_{1}$ are labelled with a circle, the nodes in $R_{2}$ are labelled
with a cross, the edges in $\gamma _{1}$ and $\gamma _{2}$ have a solid line
while those in $\gamma _{12}$ are dashed. Here we have $3+3=6$ internal sums
and six external ones.

\FRAME{dtbpFUX}{3in}{2.0003in}{0pt}{\Qcb{Figure II}}{}{Plot}{\special%
{language "Scientific Word";type "MAPLEPLOT";width 3in;height 2.0003in;depth
0pt;display "USEDEF";plot_snapshots TRUE;mustRecompute FALSE;lastEngine
"Maple";xmin "-5";xmax "5";xviewmin "-1.02083471381";xviewmax
"4.0227738179362";yviewmin "-0.9608347138";yviewmax
"0.961573817426001";plottype 4;constrained TRUE;numpoints 49;plotstyle
"patchnogrid";axesstyle "none";xis \TEXUX{x};var1name \TEXUX{$x$};function
\TEXUX{\EQN{6}{1}{}{}{\RD{\CELL{\left[ \cos \frac{\pi }{8},\sin \frac{\pi
}{8},\cos \frac{3\pi }{8},\sin \frac{3\pi }{8},\cos \frac{5\pi }{8},\sin
\frac{5\pi }{8},\cos \frac{7\pi }{8},\sin \frac{7\pi }{8},\cos \frac{9\pi
}{8},\sin \frac{9\pi }{8},\cos \frac{11\pi }{8},\sin \frac{11\pi }{8},\cos
\frac{13\pi }{8},\sin \frac{13\pi }{8},\cos \frac{15\pi }{8},\sin
\frac{15\pi }{8},\cos \frac{17\pi }{8},\sin \frac{17\pi }{8},\right]
}}{1}{}{}{}}};linecolor "black";linestyle 1;pointplot TRUE;pointstyle
"box";linethickness 1;lineAttributes "Solid";var1range
"-5,5";num-x-gridlines 49;curveColor "[flat::RGB:0000000000]";curveStyle
"Point";function \TEXUX{$\left[ \cos \frac{\pi }{8},\sin \frac{\pi }{8},\cos
\frac{9\pi }{8},\sin \frac{9\pi }{8}\right] $};linecolor "black";linestyle
1;pointstyle "point";linethickness 1;lineAttributes "Solid";var1range
"-5,5";num-x-gridlines 49;curveColor "[flat::RGB:0000000000]";curveStyle
"Line";function \TEXUX{$\left[ \cos \frac{3\pi }{8},\sin \frac{3\pi
}{8},\cos \frac{11\pi }{8},\sin \frac{11\pi }{8}\right] $};linecolor
"black";linestyle 1;pointstyle "point";linethickness 1;lineAttributes
"Solid";var1range "-5,5";num-x-gridlines 49;curveColor
"[flat::RGB:0000000000]";curveStyle "Line";function \TEXUX{$\left[ \cos
\frac{5\pi }{8},\sin \frac{5\pi }{8},\cos \frac{13\pi }{8},\sin \frac{13\pi
}{8}\right] $};linecolor "black";linestyle 1;pointstyle
"point";linethickness 1;lineAttributes "Solid";var1range
"-5,5";num-x-gridlines 49;curveColor "[flat::RGB:0000000000]";curveStyle
"Line";function \TEXUX{$\left[ \cos \frac{7\pi }{8},\sin \frac{7\pi
}{8},\cos \frac{15\pi }{8},\sin \frac{15\pi }{8},\right] $};linecolor
"black";linestyle 1;pointstyle "point";linethickness 1;lineAttributes
"Solid";var1range "-5,5";num-x-gridlines 49;curveColor
"[flat::RGB:0000000000]";curveStyle "Line";function \TEXUX{$\left[ \cos
\frac{\pi }{8},\sin \frac{\pi }{8},\cos \frac{3\pi }{8},\sin \frac{3\pi
}{8},\cos \frac{5\pi }{8},\sin \frac{5\pi }{8},\cos \frac{7\pi }{8},\sin
\frac{7\pi }{8},\cos \frac{9\pi }{8},\sin \frac{9\pi }{8},\cos \frac{11\pi
}{8},\sin \frac{11\pi }{8},\cos \frac{13\pi }{8},\sin \frac{13\pi }{8},\cos
\frac{15\pi }{8},\sin \frac{15\pi }{8},\right] $};linecolor
"black";linestyle 1;pointstyle "point";linethickness 1;lineAttributes
"Solid";var1range "-5,5";num-x-gridlines 49;curveColor
"[flat::RGB:0000000000]";curveStyle "Line";function \TEXUX{$\left[ \cos
\frac{\pi }{8},\sin \frac{\pi }{8},\cos \frac{15\pi }{8},\sin \frac{15\pi
}{8},\right] $};linecolor "black";linestyle 1;pointstyle
"point";linethickness 1;lineAttributes "Solid";var1range
"-5,5";num-x-gridlines 49;curveColor "[flat::RGB:0000000000]";curveStyle
"Line";function \TEXUX{$\left[ 3+\cos \frac{9\pi }{8},\sin \frac{9\pi
}{8},3+\cos \frac{11\pi }{8},\sin \frac{11\pi }{8},3+\cos \frac{13\pi
}{8},\sin \frac{13\pi }{8},3+\cos \frac{15\pi }{8},\sin \frac{15\pi
}{8},\right] $};linecolor "black";linestyle 1;pointstyle
"point";linethickness 1;lineAttributes "Solid";var1range
"-5,5";num-x-gridlines 49;curveColor "[flat::RGB:0000000000]";curveStyle
"Line";function \TEXUX{$\left[ 3+\cos \frac{9\pi }{8},\sin \frac{9\pi
}{8},3+\cos \frac{11\pi }{8},\sin \frac{11\pi }{8},3+\cos \frac{13\pi
}{8},\sin \frac{13\pi }{8},3+\cos \frac{15\pi }{8},\sin \frac{15\pi
}{8}\right] $};linecolor "black";linestyle 1;pointplot TRUE;pointstyle
"cross";linethickness 1;lineAttributes "Solid";var1range
"-5,5";num-x-gridlines 49;curveColor "[flat::RGB:0000000000]";curveStyle
"Point";function \TEXUX{$\left[ 3+\cos \frac{\pi }{8},\sin \frac{\pi
}{8},3+\cos \frac{3\pi }{8},\sin \frac{3\pi }{8},3+\cos \frac{5\pi }{8},\sin
\frac{5\pi }{8},3+\cos \frac{7\pi }{8},\sin \frac{7\pi }{8}\right]
$};linecolor "black";linestyle 1;pointplot TRUE;pointstyle
"circle";linethickness 1;lineAttributes "Solid";var1range
"-5,5";num-x-gridlines 49;curveColor "[flat::RGB:0000000000]";curveStyle
"Point";function \TEXUX{$\left[ 3+\cos \frac{\pi }{8},\sin \frac{\pi
}{8},3+\cos \frac{9\pi }{8},\sin \frac{9\pi }{8},\right] $};linecolor
"black";linestyle 2;pointstyle "box";linethickness 1;lineAttributes
"Dash";var1range "-5,5";num-x-gridlines 49;curveColor
"[flat::RGB:0000000000]";curveStyle "Line";function \TEXUX{$\left[ 3+\cos
\frac{5\pi }{8},\sin \frac{5\pi }{8},3+\cos \frac{13\pi }{8},\sin
\frac{13\pi }{8},\right] $};linecolor "black";linestyle 2;pointstyle
"box";linethickness 1;lineAttributes "Dash";var1range "-5,5";num-x-gridlines
49;curveColor "[flat::RGB:0000000000]";curveStyle "Line";function
\TEXUX{$\left[ 3+\cos \frac{3\pi }{8},\sin \frac{3\pi }{8},3+\cos
\frac{11\pi }{8},\sin \frac{11\pi }{8},\right] $};linecolor
"black";linestyle 2;pointstyle "box";linethickness 1;lineAttributes
"Dash";var1range "-5,5";num-x-gridlines 49;curveColor
"[flat::RGB:0000000000]";curveStyle "Line";function \TEXUX{$\left[ 3+\cos
\frac{7\pi }{8},\sin \frac{7\pi }{8},3+\cos \frac{15\pi }{8},\sin
\frac{15\pi }{8},\right] $};linecolor "black";linestyle 2;pointstyle
"box";linethickness 1;lineAttributes "Dash";var1range "-5,5";num-x-gridlines
49;curveColor "[flat::RGB:0000000000]";curveStyle "Line";function
\TEXUX{$\left[ 3+\cos \frac{\pi }{8},\sin \frac{\pi }{8},3+\cos \frac{15\pi
}{8},\sin \frac{15\pi }{8},\right] $};linecolor "black";linestyle
2;pointstyle "point";linethickness 1;lineAttributes "Dash";var1range
"-5,5";num-x-gridlines 49;curveColor "[flat::RGB:0000000000]";curveStyle
"Line";function \TEXUX{$\left[ 3+\cos \frac{7\pi }{8},\sin \frac{7\pi
}{8},3+\cos \frac{9\pi }{8},\sin \frac{9\pi }{8},\right] $};linecolor
"black";linestyle 2;pointstyle "point";linethickness 1;lineAttributes
"Dash";var1range "-5,5";num-x-gridlines 49;curveColor
"[flat::RGB:0000000000]";curveStyle "Line";function \TEXUX{$\left[ 3+\cos
\frac{\pi }{8},\sin \frac{\pi }{8},3+\cos \frac{3\pi }{8},\sin \frac{3\pi
}{8},3+\cos \frac{5\pi }{8},\sin \frac{5\pi }{8},3+\cos \frac{7\pi }{8},\sin
\frac{7\pi }{8},\right] $};linecolor "black";linestyle 1;pointstyle
"point";linethickness 1;lineAttributes "Solid";var1range
"-5,5";num-x-gridlines 49;curveColor "[flat::RGB:0000000000]";curveStyle
"Line";valid_file "T";tempfilename 'IIDRFK02.wmf';tempfile-properties "XPR";}%
}Assume now that $\gamma _{i}$ does not include any loop, for $i=1,...,g.$
Our point will be to show that%
\begin{equation}
\left| D(\gamma )\right| \leq \prod_{i=1}^{g}\left\| X_{R_{i};\gamma
_{i}}\right\| \text{ }  \label{resa}
\end{equation}%
where $\left\| .\right\| $ denotes Euclidean norm, and%
\begin{equation}
\left\| X_{R_{i};\gamma _{i}}\right\| \leq \prod_{\left\{ (r,k),.\right\}
\in \gamma _{i}}(2\ell _{rk}+1)^{-1/2}\leq (2\min_{\left\{ (r,k),.\right\}
\in \gamma _{i}}\ell _{rk}+1)^{(\#(R_{i})-1)/2};  \label{resb}
\end{equation}%
note that if $\gamma _{i}$ does not include any loop the number of edges it
contains must be identically equal to $\#(R_{i})-1,$ where $\#(.)$ denotes
the cardinality of a set. Let us consider (\ref{resa}) first. It is clear
that we can choose new indexes such that $X_{R_{i};\gamma _{i}}=:X^{(i)}$ is
a vector with elements 
\begin{equation*}
X^{(i)}=\left\{ X_{m_{i1},...,m_{iv_{i}}}^{(i)}\text{ , }-\ell _{ij}\leq
m_{ij}\leq \ell _{ij}\text{ , }j=1,...v_{i}\right\} \text{ },\text{ }%
i=1,...,g\text{ ;}
\end{equation*}%
we write $\widetilde{T}$ for this new set of indexes, namely%
\begin{equation*}
\widetilde{T}=\left\{ 
\begin{array}{cccc}
(1,1) & ... & .... & (1,v_{1}) \\ 
... & ... & ... & ... \\ 
(g,1) & ... & ... & (g,v_{g})%
\end{array}%
\right\} \text{ ;}
\end{equation*}%
here $g$ can be viewed as the number of trees and $v_{i}$ as the number of
vertexes which are in a given tree $i$. Clearly the vectors $%
X^{(1)},...,X^{(g)}$ have dimensions $\#(X^{(i)})=\prod_{j=1}^{v_{i}}(2\ell
_{ij}+1)$. The following Lemma can be viewed as an extension of the
Cauchy-Schwartz inequality.

\ \newline
\textbf{Lemma 3.1 }\emph{(Generalized Cauchy-Schwartz inequality) }Let $%
\widetilde{\gamma }$ be a partition of $\widetilde{T}$ with no flat edges. 
\begin{equation*}
\left\{ \prod_{i=1}^{g}\prod_{j=1}^{v_{i}}\sum_{m_{ij}=-\ell _{ij}}^{\ell
_{ij}}\right\} \left\{ \prod_{i=1}^{g}\left|
X_{m_{i1}...m_{i}v_{i}}^{(i)}\right| \right\} \left| \delta (\widetilde{%
\gamma })\right| \leq \prod_{i=1}^{g}\left\| X^{(i)}\right\| 
\end{equation*}%
and%
\begin{equation*}
\left\{ \prod_{i=1}^{g}\prod_{j=1}^{v_{i}}\sum_{m_{ij}=-\ell _{ij}}^{\ell
_{ij}}\right\} =\sum_{m_{11}=-\ell _{11}}^{\ell
_{11}}...\sum_{m_{gv_{g}}=-\ell _{gv_{g}}}^{\ell _{gv_{g}}}.
\end{equation*}%
\textbf{Proof }The result follows from the iterated application of the
Cauchy-Schwartz inequality; we shall argue by induction. Without loss of
generality, we can assume that the diagram $\gamma $ is connected (if it is
not, argue separately for the connected components). It is trivial to show
that the result holds for $g=2,$ indeed in that case it just the standard
Cauchy-Schwartz result. Let us now show that if the result holds for the
product of $g-1\geq 2$ components, it must hold for $g$ components as well.
Recall we consider diagrams with no flat edges, so the indexes cannot match
on the same vector $X^{(i)}.$ Relabel terms so that there exist (at least) a
link (that is, a common index $m$) between the first two vectors $%
X^{(1)},X^{(2)}$. Without loss of generality we can order terms in such a
way that the matching is internal for the first $v^{\ast }$ indexes and
external (that is, with the remaining nodes $(3,4,...,g))$ for $m_{1j}:$ $%
j=v^{\ast }+1,...,v_{1}$ and $m_{2j}:$ $j=v^{\ast }+1,...,v_{2}$. We write
also $\widetilde{\gamma }=\widetilde{\gamma }_{12}\cup \widetilde{\gamma }_{%
\overline{12}},$ where $\widetilde{\gamma }_{12}$ is the set of edges
linking node 1 to node 2 and $\widetilde{\gamma }_{\overline{12}}=\widetilde{%
\gamma }\backslash \widetilde{\gamma }_{12}$. We have 
\begin{eqnarray*}
&&\left\{ \prod_{i=1}^{g}\prod_{j=1}^{v_{i}}\sum_{m_{ij}=-\ell _{ij}}^{\ell
_{ij}}\right\} \left\{ \prod_{i=1}^{g}\left|
X_{m_{i1}...m_{i}v_{i}}^{(i)}\right| \right\} \left| \delta (\widetilde{%
\gamma })\right|  \\
&\leq &\left\{ \sum_{m_{1,v^{\ast }+1}=-\ell _{1v^{\ast }}}^{\ell _{1v^{\ast
}}}...\sum_{m_{2v_{2}}=-\ell _{2v_{2}}}^{\ell _{2v_{2}}}\right\} \left[ %
\left[ \left\{ \sum_{m_{11}=-\ell _{11}}^{\ell _{11}}...\sum_{m_{1v^{\ast
}}=-\ell _{1v^{\ast }}}^{\ell _{1v^{\ast }}}\right\} \left| \left\{
X_{m_{11}...m_{1}v_{1}}^{(1)}X_{m_{11}...m_{2}v_{2}}^{(2)}\right\} \right|
\left| \delta (\widetilde{\gamma }_{12})\right| \right] \right.  \\
&&\times \left. \left\{ \prod_{i=3}^{g}\prod_{j=1}^{v_{i}}\sum_{m_{ij}=-\ell
_{ij}}^{\ell _{ij}}\right\} \left\{ \prod_{i=3}^{g}\left|
X_{m_{i1}...m_{i}v_{i}}^{(i)}\right| \right\} \left| \delta (\widetilde{%
\gamma }_{\overline{12}})\right| \right] 
\end{eqnarray*}%
\begin{eqnarray}
&\leq &\left\{ \sum_{m_{1,v^{\ast }+1}=-\ell _{1v^{\ast }}}^{\ell _{1v^{\ast
}}}...\sum_{m_{2v_{2}}=-\ell _{2v_{2}}}^{\ell _{2v_{2}}}\right\} \left\{ 
\left[ \left\{ \sum_{m_{11}=-\ell _{11}}^{\ell _{11}}...\sum_{m_{1v^{\ast
}}=-\ell _{1v^{\ast }}}^{\ell _{1v^{\ast }}}\right\} \left\{
X_{m_{11}...m_{1}v_{1}}^{(1)}\right\} ^{2}\right] ^{1/2}\right.   \notag \\
&&\times \left[ \left\{ \sum_{m_{21}=-\ell _{21}}^{\ell
_{21}}...\sum_{m_{2v^{\ast }}=-\ell _{2v^{\ast }}}^{\ell _{2v^{\ast
}}}\right\} \left\{ X_{m_{21}...m_{2}v_{2}}^{(2)}\right\} ^{2}\right] ^{1/2}
\notag \\
&&\times \left. \left[ \left\{
\prod_{i=3}^{g}\prod_{j=1}^{v_{i}}\sum_{m_{ij}=-\ell _{ij}}^{\ell
_{ij}}\right\} \left\{ \prod_{i=3}^{g}\left|
X_{m_{i1}...m_{i}v_{i}}^{(i)}\right| \right\} \right] \left| \delta (%
\widetilde{\gamma }_{\overline{12}})\right| \right\} \text{ .}
\label{aboqua}
\end{eqnarray}%
\ The last step follows again by standard Cauchy--Schwartz inequality. Now
define 
\begin{eqnarray*}
X_{m_{1,v^{\ast }+1}...m_{2}v_{2}}^{(1;2)} &:&=\left[ \left\{
\sum_{m_{11}=-\ell _{11}}^{\ell _{11}}...\sum_{m_{1v^{\ast }}=-\ell
_{1v^{\ast }}}^{\ell _{1v^{\ast }}}\right\} \left\{
X_{m_{11}...m_{1}v_{1}}^{(1)}\right\} ^{2}\right] ^{1/2} \\
&&\times \left[ \left\{ \sum_{m_{21}=-\ell _{21}}^{\ell
_{21}}...\sum_{m_{2v^{\ast }}=-\ell _{2v^{\ast }}}^{\ell _{2v^{\ast
}}}\right\} \left\{ X_{m_{21}...m_{2}v_{2}}^{(2)}\right\} ^{2}\right] ^{1/2}
\end{eqnarray*}%
so that (\ref{aboqua}) becomes 
\begin{equation*}
\sum_{m_{1,v^{\ast }+1}=-\ell _{1v^{\ast }}}^{\ell _{1v^{\ast
}}}...\sum_{m_{2v_{2}}=-\ell _{2v_{2}}}^{\ell _{2v_{2}}}\sum_{m_{31}=-\ell
_{31}}^{\ell _{31}}...\sum_{m_{gv_{g}}=-\ell _{gv_{g}}}^{\ell
_{qv_{g}}}\left\{ X_{m_{1,v^{\ast }+1}...m_{2}v_{2}}^{(1;2)}\times
\prod_{i=3}^{g}\left| X_{m_{i1}...m_{i}v_{i}}^{(i)}\right| \right\} \left|
\delta (\widetilde{\gamma }_{\overline{12}})\right| 
\end{equation*}

\begin{equation}
\leq \left\| X_{m_{1,v^{\ast }+1}...m_{2}v_{2}}^{(1;2)}\right\|
\prod_{i=3}^{g}\left\| X^{(i)}\right\|   \label{ind}
\end{equation}%
by the inductive step (indeed within the curly brackets we have the product
of $g-2+1=g-1$ components). Now notice that\ 
\begin{equation*}
\left\| X_{m_{1,v^{\ast }+1}...m_{2}v_{2}}^{(1;2)}\right\| =\left\{
\sum_{m_{1,v^{\ast }+1}=-\ell _{1v^{\ast }}}^{\ell _{1v^{\ast
}}}...\sum_{m_{2v_{2}}=-\ell _{2v_{2}}}^{\ell _{2v_{2}}}\left(
X_{m_{1,v^{\ast }+1}...m_{2}v_{2}}^{(1;2)}\right) ^{2}\right\} ^{1/2}
\end{equation*}%
\begin{eqnarray*}
&=&\left\{ \sum_{m_{1,v^{\ast }+1}=-\ell _{1v^{\ast }}}^{\ell _{1v^{\ast
}}}...\sum_{m_{2v_{2}}=-\ell _{2v_{2}}}^{\ell _{2v_{2}}}\left\{
\sum_{m_{11}=-\ell _{11}}^{\ell _{11}}...\sum_{m_{1v^{\ast }}=-\ell
_{1v^{\ast }}}^{\ell _{1v^{\ast }}}\right\} \left\{
X_{m_{11}...m_{1}v_{1}}^{(1)}\right\} ^{2}\right.  \\
&&\left. \times \left\{ \sum_{m_{21}=-\ell _{21}}^{\ell
_{21}}...\sum_{m_{2v^{\ast }}=-\ell _{2v^{\ast }}}^{\ell _{2v^{\ast
}}}\right\} \left\{ X_{m_{21}...m_{2}v_{2}}^{(2)}\right\} ^{2}\right\} ^{1/2}
\end{eqnarray*}%
\begin{eqnarray*}
&=&\left\{ \sum_{m_{11}=-\ell _{11}}^{\ell _{11}}...\sum_{m_{1v_{1}}=-\ell
_{1v_{1}}}^{\ell _{1v_{1}}}\left\{ X_{m_{11}...m_{1}v_{1}}^{(1)}\right\}
^{2}\right\} ^{1/2} \\
&&\times \left\{ \sum_{m_{21}=-\ell _{21}}^{\ell
_{21}}...\sum_{m_{2v_{2}}=-\ell _{2v_{2}}}^{\ell _{2v_{2}}}\left\{
X_{m_{21}...m_{2}v_{2}}^{(2)}\right\} ^{2}\right\} ^{1/2} \\
&=&\left\| X^{(1)}\right\| \times \left\| X^{(2)}\right\| \text{ ,}
\end{eqnarray*}%
and thus by substitution into (\ref{ind}) the proof is completed.

\ \hfill$\square $

\ \newline
\textbf{Remark 3.2} We provide an example to make the statement of Lemma 3.1
more transparent. Take%
\begin{equation*}
X_{m_{11}m_{12}}^{(1)}\text{ , }X_{m_{21}m_{22}m_{23}}^{(2)}\text{ , }%
X_{m_{31}m_{32}m_{33}}^{(3)}\text{ , }-l_{i}\leq m_{ij}\leq l_{i}\text{ , }%
i,j=1,2,3\text{ ;}
\end{equation*}%
here%
\begin{equation*}
\widetilde{T}=\left\{ 
\begin{array}{ccc}
(1,1) & (1,2) &  \\ 
(2,1) & (2,2) & (2,3) \\ 
(3,1) & (3,2) & (3,3)%
\end{array}%
\right\} \text{ }.
\end{equation*}%
Now take for instance%
\begin{equation*}
\widetilde{\gamma }=\left[ \left\{ (1,1),(2,1)\right\} ,\left\{
(1,2),(3,1)\right\} ,\left\{ (2,2),(3,2)\right\} ,\left\{
(2,3),(3,3)\right\} \right] \text{ ,}
\end{equation*}%
and consider the sum%
\begin{eqnarray*}
&&\left\{ \prod_{i=1}^{3}\prod_{j=1}^{q_{i}}\sum_{m_{ij}}\right\}
|X_{m_{11}m_{12}}^{(1)}X_{m_{21}m_{22}m_{23}}^{(2)}X_{m_{31}m_{32}m_{33}}^{(3)}||\delta (%
\widetilde{\gamma })| \\
&=&\sum_{m_{11}}\sum_{m_{12}}\sum_{m_{22}}%
\sum_{m_{23}}|X_{m_{11}m_{12}}^{(1)}X_{m_{11}m_{22}m_{23}}^{(2)}X_{m_{12}m_{22}m_{23}}^{(3)}|
\\
&\leq &\left\| X^{(1)}\right\| \left\| X^{(2)}\right\| \left\|
X^{(3)}\right\| ,
\end{eqnarray*}%
by Lemma 3.1, where%
\begin{eqnarray*}
\left\| X^{(1)}\right\| &=&\sqrt{%
\sum_{m_{11}m_{12}}(X_{m_{11}m_{12}}^{(1)})^{2}},\text{ }\left\|
X^{(2)}\right\| =\sqrt{%
\sum_{m_{21}m_{22}m_{23}}(X_{m_{21}m_{22}m_{23}}^{(2)})^{2}}, \\
\left\| X^{(3)}\right\| &=&\sqrt{%
\sum_{m_{31}m_{32}m_{33}}(X_{m_{31}m_{32}m_{33}}^{(2)})^{2}}.
\end{eqnarray*}

\ 

It is not difficult to see that (\ref{resa}) is an immediate consequence of
Lemma 3.1; in particular, note that $\widetilde{\gamma }$ can be viewed as
the diagram which is obtained from $\gamma $ by identifying all nodes that
belong to the same set $R_{i},$ for $i=1,...,g.$ Now let us consider (\ref%
{resb}). The following result uses the previous inequality to bound the
components of a connected graph. Without loss of generality, we re-order
terms so that for all $r,$ we have $\ell _{r1}\leq \ell _{r2}\leq \ell _{r3}$
and%
\begin{equation*}
\ell _{11}\leq \ell _{21}\leq ...\leq \ell _{R1}\text{ .}
\end{equation*}%
Note that the same inequality need not be satisfied for the sequences $%
\left\{ \ell _{rk}\right\} _{r=1,...,R}$, $k=2,3.$

\ \newline
\textbf{Lemma 3.2 }

a) Every connected graph $\gamma \in \Gamma $ with no loops of orders 1,2 or
3 can be partitioned as%
\begin{equation*}
\gamma =(\gamma _{1}\cup \gamma _{2}\cup ...\cup \gamma _{g}\cup \gamma
_{12}\cup ...\cup \gamma _{g-1,g})
\end{equation*}%
where $\gamma _{i}$ has no loops of any order and is such that $%
\#(R_{i})\geq 2$ (in other words, $\gamma $ can be broken into into binary
trees with at least two nodes)

b) For all $\gamma _{i},$ $i=1,...,g$ we have

\begin{equation*}
\left\| X_{R_{i};\gamma _{i}}\right\| \leq \prod_{\left\{ (r,k),.\right\}
\in \gamma _{i}}(2\ell _{rk}+1)^{-1/2}\leq (2\min_{\left\{ (r,k),.\right\}
\in \gamma _{i}}\ell _{rk}+1)^{(\#(R_{i})-1)/2}
\end{equation*}%
that is, every tree with $p$ nodes corresponds to summands which are $%
O((2\ell _{11}+1)^{-(p-1)/2})$, where $p$ is the number of nodes in the
trees.

c) For all connected graphs $\gamma $ with $p$ nodes, we have 
\begin{equation}
\left| D[\gamma ;\mathcal{L}]\right| \leq \prod_{r=1}^{p/4}(2\ell
_{r1}+1)^{-1}.  \label{bound}
\end{equation}%
\newline
\textbf{Proof}

a) We drop edges till we reach the point where there are only binary trees
or isolated points. Each of these isolated points can be connected to either
another isolated point, in which case we simply have a tree with two nodes,
or to a binary tree. The graph which is obtained by linking this point to
the tree is itself a tree if its paths have at most two edges: recall there
are no loops of order 2 or 3. On the other hand, if the graph has a path
which covers four nodes, then we delete one edge and obtain two trees with
two nodes. These procedures can be iterated until no isolated point remains.

b) We recall the identities (see VMK, chapter 8)%
\begin{eqnarray*}
\left( 
\begin{array}{ccc}
\ell _{1} & \ell _{2} & \ell _{3} \\ 
m_{1} & m_{2} & -m_{3}%
\end{array}%
\right) &=&(-1)^{\ell _{3}+m_{3}}\frac{1}{\sqrt{2\ell _{3}+1}}C_{\ell
_{1},-m_{1},\ell _{2},-m_{2}}^{\ell _{3},m_{3}} \\
C_{\ell _{1},m_{1},\ell _{2},m_{2}}^{\ell _{3},m_{3}} &=&(-1)^{\ell
_{1}-\ell _{2}+m_{3}}\sqrt{2\ell _{3}+1}\left( 
\begin{array}{ccc}
\ell _{1} & \ell _{2} & \ell _{3} \\ 
m_{1} & m_{2} & -m_{3}%
\end{array}%
\right) \text{ ,}
\end{eqnarray*}%
where the coefficients $C_{\ell _{1},m_{1},\ell _{2},m_{2}}^{\ell
_{3},m_{3}} $ (\emph{Clebsch-Gordan coefficients}) are the elements of a
unitary matrix which implements the change of basis from a tensor product to
a direct sum representation for a space spanned by spherical harmonics. More
precisely, denote by $\left\{ Y_{\ell _{1}}\otimes Y_{\ell _{2}}\right\}
_{m_{1},m_{2}}$ the elements of a basis for the tensor product $Y_{\ell
_{1}}\otimes Y_{\ell _{2}};$ here, $Y_{\ell }$ is the vector space generated
by $Y_{\ell m},$ $m=-\ell ,...,\ell :$ see for instance Vilenkin and Klimyk
(1991), Chapter 8 for a full discussion of tensor products and their
properties. For our purposes, it suffices to recall the identity

\begin{equation*}
\left\{ Y_{\ell _{1}}\otimes Y_{\ell _{2}}\right\} _{m_{1},m_{2}}=\sum_{\ell
=|\ell _{2}-\ell _{1}|}^{\ell _{2}+\ell _{1}}\sum_{m=-\ell }^{\ell }C_{\ell
_{1},m_{1},\ell _{2},m_{2}}^{\ell ,m}Y_{\ell m}\text{ .}
\end{equation*}%
More compactly, we might refer to the $(2\ell _{1}+1)(2\ell _{2}+1)\times
(2\ell _{1}+1)(2\ell _{2}+1)$ matrix $\mathcal{C}$, whose elements $\left\{
C_{\ell _{1},m_{1},\ell _{2},m_{2}}^{\ell ,m}\right\} $ are indexed by $%
m_{1},m_{2}$ over the rows and $\ell ,m$ over the columns. The matrix $%
\mathcal{C}$ is unitary: it transforms a basis the orthonormality
relationships read%
\begin{eqnarray*}
\sum_{m_{1},m_{2}}C_{\ell _{1},m_{1},\ell _{2},m_{2}}^{\ell ,m}C_{\ell
_{1},m_{1},\ell _{2},m_{2}}^{\ell ^{\prime },m^{\prime }} &=&\delta _{\ell
}^{\ell \prime }\delta _{m}^{m\prime }, \\
\sum_{\ell ,m}C_{\ell _{1},m_{1},\ell _{2},m_{2}}^{\ell ,m}C_{\ell
_{1},m_{1}^{\prime },\ell _{2},m_{2}^{\prime }}^{\ell ,m} &=&\delta
_{m_{1}}^{m_{1}^{\prime }}\delta _{m_{2}}^{m_{2}^{\prime }}\text{ .}
\end{eqnarray*}%
Now the argument can be iterated to higher-order tensor products, to obtain%
\begin{eqnarray*}
\left\{ Y_{\ell _{1}}\otimes Y_{\ell _{2}}\otimes Y_{\ell _{3}}\right\}
_{m_{1},m_{2},m_{3}} &=&\sum_{\ell =|\ell _{2}-\ell _{1}|}^{\ell _{2}+\ell
_{1}}\sum_{m=-\ell }^{\ell }C_{\ell _{1},m_{1},\ell _{2},m_{2}}^{\ell
,m}\left\{ Y_{\ell }\otimes Y_{\ell _{3}}\right\} _{m,m_{3}} \\
&=&\sum_{\ell =|\ell _{2}-\ell _{1}|}^{\ell _{2}+\ell _{1}}\sum_{\ell
^{\prime }=|\ell _{3}-\ell |}^{\ell _{3}+\ell }\sum_{m=-\ell }^{\ell
}\sum_{m^{\prime }=-\ell ^{\prime }}^{\ell \prime }C_{\ell _{1},m_{1},\ell
_{2},m_{2}}^{\ell ,m}C_{\ell ,m,\ell _{3},m_{3}}^{\ell ^{\prime },m^{\prime
}}Y_{\ell ^{\prime }m^{\prime }}\text{ .}
\end{eqnarray*}%
The orthonormality conditions now read%
\begin{equation*}
\sum_{m_{1}m_{2}m_{3}}\left\{ \sum_{m=-\ell }^{\ell }C_{\ell _{1},m_{1},\ell
_{2},m_{2}}^{\ell ,m}C_{\ell ,m,\ell _{3},m_{3}}^{\ell ^{\prime },m^{\prime
}}\right\} ^{2}=\sum_{\ell \ell ^{\prime }m^{\prime }}\left\{ \sum_{m=-\ell
}^{\ell }C_{\ell _{1},m_{1},\ell _{2},m_{2}}^{\ell ,m}C_{\ell ,m,\ell
_{3},m_{3}}^{\ell ^{\prime },m^{\prime }}\right\} ^{2}=1\text{ .}
\end{equation*}%
More precisely, the coefficients 
\begin{equation*}
C_{\ell _{1},m_{1},\ell _{2},m_{2},\ell _{3}m_{3}}^{\ell ,\ell ^{\prime
},m}:=\sum_{m=-\ell }^{\ell }C_{\ell _{1},m_{1},\ell _{2},m_{2}}^{\ell
,m}C_{\ell ,m,\ell _{3},m_{3}}^{\ell ^{\prime },m^{\prime }}
\end{equation*}%
are the elements of a unitary matrix with $(2\ell _{1}+1)(2\ell
_{2}+1)(2\ell _{3}+1)$ rows and columns; the columns are indexed by $\ell
,\ell ^{\prime },m$ and they are%
\begin{equation*}
\sum_{\ell =|\ell _{2}-\ell _{1}|}^{\ell _{2}+\ell _{1}}\sum_{\ell ^{\prime
}=|\ell _{3}-\ell |}^{\ell _{3}+\ell }(2\ell ^{\prime }+1)=(2\ell
_{3}+1)\sum_{\ell =|\ell _{2}-\ell _{1}|}^{\ell _{2}+\ell _{1}}(2\ell
+1)=(2\ell _{1}+1)(2\ell _{2}+1)(2\ell _{3}+1)\text{ .}
\end{equation*}%
In general, we have the coefficients%
\begin{equation*}
C_{\ell _{1},m_{1};...;\ell _{p}m_{p}}^{\lambda _{1},\lambda
_{2},...,\lambda _{p-1};\mu }:=\sum_{\mu _{1}=-\lambda _{1}}^{\lambda
_{1}}...\sum_{\mu _{p-2}=-\lambda _{p-2}}^{\lambda _{p-2}}C_{\ell
_{1},m_{1},\ell _{2},m_{2}}^{\lambda _{1},\mu _{1}}\times C_{\lambda
_{1},\mu _{1};\ell _{3},m_{3}}^{\lambda _{2},\mu _{2}}\times ...\times
C_{\lambda _{p-2},\mu _{p-2};\ell _{p},m_{p}}^{\lambda _{p-1},\mu }
\end{equation*}%
such that 
\begin{equation*}
\left\{ Y_{\ell _{1}}\otimes ...\otimes Y_{\ell _{p}}\right\}
_{m_{1},...;m_{p}}=\sum_{\lambda _{1}}...\sum_{\lambda _{p-2}}\sum_{\lambda
_{p-1}}\sum_{\mu =-\lambda _{p-1}}^{\lambda _{p-1}}C_{\ell
_{1},m_{1};...;\ell _{p}m_{p}}^{\lambda _{1},\lambda _{2},...,\lambda
_{p-1};\mu }Y_{\lambda _{p-1},\mu }\text{ .}
\end{equation*}%
whence the orthonormality conditions yield 
\begin{equation*}
\sum_{m_{1},...m_{p}}\left\{ C_{\ell _{1},m_{1};...;\ell _{p}m_{p}}^{\lambda
_{1},\lambda _{2},...,\lambda _{p-1};\mu }\right\} ^{2}=\sum_{\lambda
_{1}}...\sum_{\lambda _{p-2}}\sum_{\lambda _{p-1}}\sum_{\mu =-\lambda
_{p-1}}^{\lambda _{p-1}}\left\{ C_{\ell _{1},m_{1};...;\ell
_{p}m_{p}}^{\lambda _{1},\lambda _{2},...,\lambda _{p-1};\mu }\right\} ^{2}=1%
\text{ .}
\end{equation*}%
It follows that any cross-product of $p$ Wigner's $3j$ coefficients
corresponding to a binary tree can be bounded by a term like%
\begin{eqnarray*}
&&\sum_{m_{1},...m_{p};m}\left| \sum_{\mu _{1}=-\lambda _{1}}^{\lambda
_{1}}...\sum_{\mu _{p-1}=-\lambda _{p-1}}^{\lambda _{p-1}}\left( 
\begin{tabular}{lll}
$\ell _{1}$ & $\ell _{2}$ & $\lambda _{1}$ \\ 
$m_{1}$ & $m_{2}$ & $\mu _{1}$%
\end{tabular}%
\right) \left( 
\begin{tabular}{lll}
$\lambda _{1}$ & $\ell _{3}$ & $\lambda _{2}$ \\ 
$\mu _{1}$ & $m_{3}$ & $\mu _{2}$%
\end{tabular}%
\right) ....\left( 
\begin{tabular}{lll}
$\lambda _{p-1}$ & $\ell _{p+1}$ & $\ell $ \\ 
$\mu _{p-1}$ & $m_{p+1}$ & $m$%
\end{tabular}%
\right) \right| ^{2} \\
&=&\left\{ (2\ell +1)\prod_{j=1}^{p-1}(2\lambda _{j}+1)\right\} ^{-1} \\
&&\times \sum_{m_{1},...m_{p};m}\left| \sum_{\mu _{1}=-\lambda
_{1}}^{\lambda _{1}}...\sum_{\mu _{p-1}=-\lambda _{p-1}}^{\lambda
_{p-1}}C_{\ell _{1},m_{1},\ell _{2},m_{2}}^{\lambda _{1},\mu _{1}}\times
C_{\lambda _{1},\mu _{1};\ell _{3},m_{3}}^{\lambda _{2},\mu _{2}}\times
...\times C_{\lambda _{p-1},\mu _{p-1};\ell _{p},m_{p}}^{\ell ,m}\right| ^{2}
\\
&=&\left\{ (2\ell +1)\prod_{j=1}^{p-1}(2\lambda _{j}+1)\right\}
^{-1}\sum_{m_{1},...m_{p};m}\left\{ C_{\ell _{1},m_{1};...;\ell
_{p}m_{p}}^{\lambda _{1},\lambda _{2},...,\lambda _{p-1},\ell ;m}\right\}
^{2}=\left\{ \prod_{j=1}^{p-1}(2\lambda _{j}+1)\right\} ^{-1}.
\end{eqnarray*}

c) Connected components with loops of order 1, 2 or 3 can be reduced as
shown in Marinucci (2005). From part b), we have 
\begin{equation}
\left| D[\gamma ;\mathcal{L}]\right| \leq \prod_{r=1}^{(p-g)/2}(2\ell
_{r1}+1)^{-1},  \label{bound2}
\end{equation}%
whereas from a) we learn that it is always possible to choose a partition
such that $g\leq p/2$; the result follows immediately.

\ \hfill$\square $

\ \newline
\textbf{Remark 3.3} An inspection of our argument (especially in c)) reveals
that the bound (\ref{bound}) can be improved in those cases where it is
possible to partition the graph $\gamma $ into a minimal number of trees.
For instance, it is known (see for instance Biederharn and Louck (1981))
that all connected graphs with up to 8 nodes can be partitioned into two
binary trees; for such cases, we hence obtain the bound 
\begin{equation}
\left| D[\gamma ;\mathcal{L}]\right| \leq \prod_{r=1}^{p/2-1}(2\ell
_{r1}+1)^{-1},\text{ }p\leq 8\text{ .}  \label{bound3}
\end{equation}%
For $p\geq 10$ it is no longer true that such a partition necessarily
exists; we leave for future research, however, to ascertain whether (\ref%
{bound2}) can be improved by a more efficient use of graphical arguments.\ 

\ 

We consider now the case where the angular power spectrum is unknown and
estimated from the data. Define 
\begin{equation}
u_{lm}:=\frac{a_{lm}}{\sqrt{C_{l}}}\text{ , }\widehat{u}_{lm}:=\frac{a_{lm}}{%
\sqrt{\widehat{C}_{l}}}\text{ , }m=0,1,...,l\text{ ;}  \label{def3}
\end{equation}%
from Marinucci (2005) we recall the following simple result.

\qquad \newline
\textbf{Lemma 3.3 }Let $l$ and $p$ be positive integers, and define%
\begin{equation}
g(l;p):=\prod_{k=1}^{p}\left\{ \frac{2l+1}{2l+2k-1}\right\} \text{ .}
\label{prevg}
\end{equation}%
For $u$ and $\widehat{u}$ defined by (\ref{def3}), we have%
\begin{equation*}
E\left\{ \underset{q_{0}\text{ times}}{\underbrace{\widehat{u}_{l0}...%
\widehat{u}_{l0}}}\underset{q_{1}\text{ times}}{\underbrace{\widehat{u}%
_{l1}...\widehat{u}_{l1}}}\underset{q_{1}^{\prime }\text{ times}}{%
\underbrace{\widehat{u}_{l1}^{\ast }...\widehat{u}_{l1}^{\ast }}}...\underset%
{q_{k}\text{ times}}{\underbrace{\widehat{u}_{lk}...\widehat{u}_{lk}}}%
\underset{q_{k}^{\prime }\text{ times}}{\underbrace{\widehat{u}_{lk}^{\ast
}...\widehat{u}_{lk}^{\ast }}}\right\}
\end{equation*}%
\begin{equation*}
=E\left\{ \underset{q_{0}\text{ times}}{\underbrace{u_{l0}...u_{l0}}}%
\underset{q_{1}\text{ times}}{\underbrace{u_{l1}...u_{l1}}}\underset{%
q_{1}^{\prime }\text{ times}}{\underbrace{u_{l1}^{\ast }...u_{l1}^{\ast }}}%
...\underset{q_{k}\text{ times}}{\underbrace{u_{lk}...u_{lk}}}\underset{%
q_{k}^{\prime }\text{ times}}{\underbrace{u_{lk}^{\ast }...u_{lk}^{\ast }}}%
\right\} \times g(l;q_{0}+q_{1}+...+q_{k}^{\prime })\text{ .}
\end{equation*}

\ 

Now partition the indexes $\left\{ \ell _{11},\ell _{12},\ell _{13},...,\ell
_{R1},\ell _{R2},\ell _{R3}\right\} $ into equivalence classes where $\ell
_{rk}$ takes the same value; we can then relabel them as $\overline{\ell }%
_{u}$, $u=1,...,c$ where $c$ is the number of such classes, i.e. the number
of different values of $\overline{\ell };$ clearly $c\leq 3R.$ For
convenience, we write%
\begin{equation*}
\left[ \left\{
(l_{11},l_{12},l_{13}),...(l_{I1},l_{I2},l_{I3});p_{1},...,p_{I}\right\} %
\right] :=\left\{ \ell _{11},\ell _{12},\ell _{13},...,\ell _{R1},\ell
_{R2},\ell _{R3}\right\} \text{ ,}
\end{equation*}%
and define 
\begin{eqnarray*}
G\left[ \left\{
(l_{11},l_{12},l_{13}),...(l_{I1},l_{I2},l_{I3});p_{1},...,p_{I}\right\} %
\right] &:&=\prod_{u=1}^{c}g(\overline{\ell }_{u};\frac{\#(u)}{2}) \\
&=&\prod_{u=1}^{c}\prod_{j=1}^{(\#(u))/2}\left\{ \frac{2\overline{\ell }%
_{u}+1}{2\overline{\ell }_{u}+2j-1}\right\} \text{ ,}
\end{eqnarray*}%
where $\#(u)$ denotes the cardinality of the class where $\overline{\ell }%
_{u}$ belongs, or in other words the number of times that each particular
value is repeated in $\left\{ \ell _{11},\ell _{12},\ell _{13},...,\ell
_{R1},\ell _{R2},\ell _{R3}\right\} $. We give some example to make the
previous definition more transparent. Consider for notational simplicity the
univariate case $I=1.$ For $l_{1}<l_{2}<l_{3}$ we have $c=3$ and%
\begin{equation*}
G\left[ \left\{
(l_{11},l_{12},l_{13}),...(l_{I1},l_{I2},l_{I3});p_{1},...,p_{I}\right\} %
\right] =\prod_{u=1}^{3}\prod_{j=1}^{R}\left\{ \frac{2l_{u}+1}{2l_{u}+2j-1}%
\right\} \text{ ;}
\end{equation*}%
for $l_{1}=l_{2}<l_{3}$ we have $c=2$ and

\begin{equation*}
G\left[ \left\{
(l_{11},l_{12},l_{13}),...(l_{I1},l_{I2},l_{I3});p_{1},...,p_{I}\right\} %
\right] =\prod_{j=1}^{2R}\left\{ \frac{2l_{1}+1}{2l_{1}+2j-1}\right\}
\prod_{j=1}^{R}\left\{ \frac{2l_{3}+1}{2l_{3}+2j-1}\right\} \text{ ;}
\end{equation*}%
for $l_{1}=l_{2}=l_{3}$ we have $c=1$ and%
\begin{equation*}
G\left[ \left\{
(l_{11},l_{12},l_{13}),...(l_{I1},l_{I2},l_{I3});p_{1},...,p_{I}\right\} %
\right] =\prod_{j=1}^{3R}\left\{ \frac{2l_{1}+1}{2l_{1}+2j-1}\right\} \text{
.}
\end{equation*}%
It is shown in Marinucci (2005) that (see (\ref{prevg}))%
\begin{equation*}
E\widehat{I}%
_{l_{1}l_{2}l_{3}}^{2p}=EI_{l_{1}l_{2}l_{3}}^{2p}g(l_{1},l_{2},l_{3};p)\text{
;}
\end{equation*}%
for any $p\in \mathbb{N}$; likewise, take $p_{1},...,p_{I}\in \mathbb{N}$:
by exactly the same argument it is immediate to obtain%
\begin{equation}
E\left\{ \prod_{i=1}^{I}\widehat{I}_{l_{i1}l_{i2}l_{i3}}^{2p_{i}}\right\}
=E\left\{ \prod_{i=1}^{I}I_{l_{i1}l_{i2}l_{i3}}^{2p_{i}}\right\} G\left[
\left\{
(l_{11},l_{12},l_{13}),...(l_{I1},l_{I2},l_{I3});p_{1},...,p_{I}\right\} %
\right] \text{ .}  \label{formest}
\end{equation}%
Notice that.%
\begin{equation*}
\lim_{l_{11}\rightarrow \infty }G\left[ \left\{
(l_{11},l_{12},l_{13}),...(l_{I1},l_{I2},l_{I3});p_{1},...,p_{I}\right\} %
\right] =1\text{ .}
\end{equation*}

\ \newline
The following result extends Theorem 3.1 to the case where the angular power
spectrum is estimated from the data.

\ \newline
\textbf{Theorem 3.2 }Assume that $(l_{i1},l_{i2},l_{i3})\neq (l_{i^{\prime
}1},l_{i^{\prime }2},l_{i^{\prime }3})$ whenever $i\neq i^{\prime }.$ There
exist an absolute constant $K_{p_{1}...p_{I}}$ such that, for $p_{i}\geq 1$
, $i=1,...,I$%
\begin{equation*}
\left| E\left\{ \prod_{i=1}^{I}\widehat{I}_{l_{i1}l_{i2}l_{i3}}^{2p_{i}}%
\right\} -\prod_{i=1}^{I}\left\{ (2p_{i}-1)!!\Delta
_{l_{i1}l_{i2}l_{i3}}^{p_{i}}\right\} G\left[ \left\{
(l_{11},l_{12},l_{13}),...(l_{I1},l_{I2},l_{I3});p_{1},...,p_{I}\right\} %
\right] \right| \leq \frac{K_{p_{1}...p_{I}}}{2l_{11}+1}\text{ .}
\end{equation*}%
\textbf{Proof} It is sufficient to combine Theorem 3.1 and (\ref{formest}).

\ \hfill$\square $

\ \newline
An immediate consequence of the Theorems 3.1 and 3.2 is the following

\ \newline
\textbf{Theorem 3.3 }\emph{(Multivariate Central Limit Theorem)}\textbf{\ }%
For any $k\in \mathbb{N},$ as $l_{11}\rightarrow \infty ,$%
\begin{equation*}
\left( \frac{I_{l_{11}l_{12}l_{13}}}{\sqrt{\Delta _{l_{11}l_{12}l_{13}}}}%
,...,\frac{I_{l_{k1}l_{k2}l_{k3}}}{\sqrt{\Delta _{l_{l1}l_{k2}l_{k3}}}}%
\right) ,\left( \frac{\widehat{I}_{l_{11}l_{12}l_{13}}}{\sqrt{\Delta
_{l_{11}l_{12}l_{13}}}},...,\frac{\widehat{I}_{l_{k1}l_{k2}l_{k3}}}{\sqrt{%
\Delta _{l_{l1}l_{k2}l_{k3}}}}\right) \rightarrow _{d}N(0,I_{k})\text{ ,}
\end{equation*}%
where $I_{k}$ denotes the $(k\times 1)$ identity matrix.\newline
\textbf{Proof} The results follow immediately from Theorems 3.1, 3.2 and the
method of moments.

\ \hfill$\square $

\ \newline
\textbf{Remark 3.4} It is immediate to see that, for some $%
K_{p_{1}...p_{I}}^{\prime },K_{p_{1}...p_{I}}^{\prime \prime }>0$ 
\begin{equation*}
\left| G\left[ \left\{
(l_{11},l_{12},l_{13}),...(l_{I1},l_{I2},l_{I3});p_{1},...,p_{I}\right\} %
\right] -1\right| \leq \frac{K_{p_{1}...p_{I}}^{\prime }}{2l_{11}+1}
\end{equation*}%
whence Theorem 3.1 can be also formulated as 
\begin{equation*}
\left| E\left\{ \prod_{i=1}^{I}\widehat{I}_{l_{i1}l_{i2}l_{i3}}^{2p_{i}}%
\right\} -\prod_{i=1}^{I}\left\{ (2p_{i}-1)!!\Delta
_{l_{i1}l_{i2}l_{i3}}^{p_{i}}\right\} \right| \leq \frac{K_{p_{1}...p_{I}}^{%
\prime \prime }}{2l_{11}+1}\text{ .}
\end{equation*}%
A careful inspection of the proofs reveals that the rates of the bounds are
exact (we shall come back to this point in the next Section); in other
words, the bispectrum converges to a Gaussian distribution with the same
rate when either the angular power spectrum is known or unknown (the
bounding constants differ, however). This result settles some questions
raised in Komatsu et al. (2002)), where the distributions of $%
I_{l_{i1}l_{i2}l_{i3}}$ and $\widehat{I}_{l_{i1}l_{i2}l_{i3}}$ where
compared by means of Monte Carlo simulations.

\ \ \newline
\textbf{Remark 3.5 }It is interesting to note that the angular bispectrum
ordinates at different multipoles are asymptotically independent, for any
triples $(l_{11},l_{12},l_{13})\neq (l_{21},l_{22},l_{23}).$

\ 

\begin{center}
\textbf{4. HIGHER ORDER RESULTS}
\end{center}

\ 

The results of the previous sections can be further developed to provide
some higher order approximation on the moments of the angular bispectrum.
For notational simplicity, in this Section we shall focus on the univariate
case, that is, we do not consider cross-moments; the multivariate
generalization is straightforward, however, and no new ideas are required.
In the statement of the Theorem to follow, we use the Wigner's $6j$ symbols,
which are defined by 
\begin{equation}
\left\{ 
\begin{array}{ccc}
a & b & e \\ 
c & d & f%
\end{array}%
\right\} :=\sum_{\alpha ,\beta ,\gamma }\sum_{\varepsilon ,\delta ,\phi
}(-1)^{e+f+\varepsilon +\phi }\left( 
\begin{array}{ccc}
a & b & e \\ 
\alpha & \beta & \varepsilon%
\end{array}%
\right) \left( 
\begin{array}{ccc}
c & d & e \\ 
\gamma & \delta & -\varepsilon%
\end{array}%
\right) \left( 
\begin{array}{ccc}
a & d & f \\ 
\alpha & \delta & -\phi%
\end{array}%
\right) \left( 
\begin{array}{ccc}
c & b & f \\ 
\gamma & \beta & \phi%
\end{array}%
\right) \text{ },  \label{6j1}
\end{equation}
see VMK, chapter 9 for a full set of properties; we simply recall here that 
\begin{equation}
\left| \left\{ 
\begin{tabular}{lll}
$a$ & $b$ & $e$ \\ 
$c$ & $d$ & $f$%
\end{tabular}%
\right\} \right| \leq \min \left( \frac{1}{\sqrt{(2a+1)(2c+1)}},\frac{1}{%
\sqrt{(2b+1)(2d+1)}},\frac{1}{\sqrt{(2e+1)(2f+1)}}\right) \text{ .}
\label{6j2}
\end{equation}

\ \newline
\textbf{Theorem 4.1 }There exist an absolute constant $K_{p}$ such that

\begin{equation}
\left| \frac{EI_{l_{1}l_{2}l_{3}}^{2p}}{\Delta _{l_{1}l_{2}l_{3}}^{p}}%
-(2p-1)!!-\left\{ \frac{p(p-1)}{6}\kappa _{4}(l_{1},l_{2},l_{3})\right\}
(2p-1)!!\right| \leq \frac{K_{p}}{(2l_{1}+1)^{2}}\text{ ,}  \label{teo4.1}
\end{equation}%
where 
\begin{equation*}
\kappa _{4}(l_{1},l_{2},l_{3}):=\frac{6}{2l_{1}+1}+\frac{6}{2l_{2}+1}+\frac{6%
}{2l_{3}+1}+6\left\{ 
\begin{tabular}{lll}
$l_{1}$ & $l_{2}$ & $l_{3}$ \\ 
$l_{1}$ & $l_{2}$ & $l_{3}$%
\end{tabular}%
\right\} \text{ , for}\ l_{1}<l_{2}<l_{3}\text{ ,}
\end{equation*}

\begin{eqnarray*}
\kappa _{4}(l_{1},l_{2},l_{3}) &:&=\frac{96}{2l_{2}+1}+\frac{24}{2l_{3}+1}%
+48\left\{ 
\begin{tabular}{lll}
$l_{1}$ & $l_{2}$ & $l_{3}$ \\ 
$l_{1}$ & $l_{2}$ & $l_{3}$%
\end{tabular}%
\right\} \text{ , for }l_{1}=l_{2}<l_{3}\text{ ,} \\
\kappa _{4}(l_{1},l_{2},l_{3}) &:&=\frac{96}{2l_{2}+1}+\frac{24}{2l_{1}+1}%
+48\left\{ 
\begin{tabular}{lll}
$l_{1}$ & $l_{2}$ & $l_{3}$ \\ 
$l_{1}$ & $l_{2}$ & $l_{3}$%
\end{tabular}%
\right\} \text{ , for }l_{1}<l_{2}=l_{3}\text{ ,}
\end{eqnarray*}

\begin{equation*}
\kappa _{4}(l,l,l):=\frac{6\times 18^{2}}{2l+1}+6^{4}\left\{ 
\begin{tabular}{lll}
$l$ & $l$ & $l$ \\ 
$l$ & $l$ & $l$%
\end{tabular}%
\right\} \text{ , for }l_{1}=l_{2}=l_{3}=l\text{ .}
\end{equation*}%
It holds that (see (\ref{6j2}))%
\begin{equation*}
\left| \kappa _{4}(l_{1},l_{2},l_{3})\right| \leq \frac{C}{2l_{1}+1}\text{ ,
some }C>0\text{ }.
\end{equation*}%
\newline
\textbf{Proof} For all $p\in \mathbb{N},$ a direct combinatorial argument
yields the following

\begin{equation*}
EI_{l_{1}l_{2}l_{3}}^{2p}=(2p-1)!!\Delta _{l_{1}l_{2}l_{3}}^{p}\text{ }%
(=:D_{1}(l_{1},l_{2},l_{3}))
\end{equation*}%
\begin{equation*}
\left. +\left( 
\begin{tabular}{l}
$2p$ \\ 
$4$%
\end{tabular}%
\right) \kappa _{4}(l_{1},l_{2},l_{3})(2p-5)!!\Delta
_{l_{1}l_{2}l_{3}}^{p-2}\right\} \text{ }(=:D_{2}(l_{1},l_{2},l_{3}))
\end{equation*}%
\begin{equation*}
\left. 
\begin{array}{c}
+\frac{1}{2}\left( 
\begin{tabular}{l}
$2p$ \\ 
$4$%
\end{tabular}%
\right) \left( 
\begin{tabular}{l}
$2p-4$ \\ 
$4$%
\end{tabular}%
\right) \kappa _{4}^{2}(l_{1},l_{2},l_{3})(2p-9)!!\Delta
_{l_{1}l_{2}l_{3}}^{p-4} \\ 
+\left( 
\begin{tabular}{l}
$2p$ \\ 
$6$%
\end{tabular}%
\right) \kappa _{6}(l_{1},l_{2},l_{3})(2p-7)!!\Delta _{l_{1}l_{2}l_{3}}^{p-3}%
\end{array}%
\right\} \text{ }(=:D_{3}(l_{1},l_{2},l_{3}))
\end{equation*}%
\begin{equation*}
\left. 
\begin{array}{c}
+\frac{1}{3!}\left( 
\begin{tabular}{l}
$2p$ \\ 
$4$%
\end{tabular}%
\right) \left( 
\begin{tabular}{l}
$2p-4$ \\ 
$4$%
\end{tabular}%
\right) \left( 
\begin{tabular}{l}
$2p-8$ \\ 
$4$%
\end{tabular}%
\right) \kappa _{4}^{3}(l_{1},l_{2},l_{3})\times (2p-13)!!\Delta
_{l_{1}l_{2}l_{3}}^{p-6} \\ 
+\left( 
\begin{tabular}{l}
$2p$ \\ 
$6$%
\end{tabular}%
\right) \left( 
\begin{tabular}{l}
$2p-6$ \\ 
$4$%
\end{tabular}%
\right) \left( 
\begin{tabular}{l}
$2p-10$ \\ 
$2$%
\end{tabular}%
\right) \kappa _{6}(l_{1},l_{2},l_{3})\kappa _{4}(l_{1},l_{2},l_{3})\times
(2p-11)!!\Delta _{l_{1}l_{2}l_{3}}^{p-5} \\ 
+\left( 
\begin{tabular}{l}
$2p$ \\ 
$8$%
\end{tabular}%
\right) \kappa _{8}(l_{1},l_{2},l_{3})(2p-9)!!\Delta _{l_{1}l_{2}l_{3}}^{p-4}%
\end{array}%
\right\} \text{ }(=:D_{4}(l_{1},l_{2},l_{3}))
\end{equation*}%
\begin{equation*}
+...+\kappa _{2p}(l_{1},l_{2},l_{3})\text{ }(=:D_{p}(l_{1},l_{2},l_{3}))%
\text{ ,}
\end{equation*}%
where the sum runs over all elements such that the factorials are
nonnegative; also, we take $(2p-k)!!=1$ when $k>2p$ and%
\begin{equation*}
\kappa _{u}(l_{1},l_{2},l_{3}):=\sum_{\gamma }D[\gamma ;l_{1},l_{2},l_{3}]%
\text{ ;}
\end{equation*}%
where the sum runs over all possible connected graphs with $u$ nodes, that
is $\kappa _{u}(l_{1},l_{2},l_{3})$ represents the expected value
corresponding to connected components with $u$ nodes$.$ The term $%
D_{1}(l_{1},l_{2},l_{3})$ corresponds to the sums over all graphs with
exactly $p$ connected components, each of them with exactly two nodes. The
term $D_{2}(l_{1},l_{2},l_{3})$ correspond to the sums over all graphs with
exactly $p-1$ connected components: the number of such graphs corresponds to
the possible ways to select $4$ nodes out of $p$ and then partition the
remaining $2p-4$ into pairs, that is.%
\begin{equation*}
\left( 
\begin{tabular}{l}
$2p$ \\ 
$4$%
\end{tabular}%
\right) \times (2p-5)!!\Delta _{l_{1}l_{2}l_{3}}^{p-2}\text{ .}
\end{equation*}%
The argument for the remaining terms is entirely analogous. Now from the
previous section we know that $\kappa _{u}(l_{1},l_{2},l_{3})=O(l_{1}^{-\min
(2,\frac{u}{4})})$ for $u\geq 6$, whence the result will follow from an
explicit evaluation of $\kappa _{4}(l_{1},l_{2},l_{3}).$.For the latter,
note that connected graphs with four nodes must include a loop of order 1, 2
or 3. Terms with a 1-loop are identically zero (Marinucci (2005), Lemma
3.1). Graphs associated to terms with 2- or 3-loops ($\gamma \in \Gamma
_{L(2)}$ or $\gamma \in \Gamma _{L(3)})$ are represented in Figure III.%
\FRAME{dtbpFUX}{3in}{2.0003in}{0pt}{\Qcb{Figure III: $\protect\gamma \in
\Gamma _{L(2)},\Gamma _{L(3)}$}}{}{Plot}{\special{language "Scientific
Word";type "MAPLEPLOT";width 3in;height 2.0003in;depth 0pt;display
"USEDEF";plot_snapshots TRUE;mustRecompute FALSE;lastEngine "Maple";xmin
"1";xmax "2";xviewmin "0.94";xviewmax "4.0612";yviewmin
"0.850000190031882";yviewmax "2.15049980982194";plottype 4;constrained
TRUE;numpoints 49;plotstyle "patch";axesstyle "none";xis \TEXUX{x};var1name
\TEXUX{$x$};function
\TEXUX{\EQN{6}{1}{}{}{\RD{%
\CELL{((x-1)^{1.5}(2-x)^{1.5}+1)}}{1}{}{}{}}};linecolor "black";linestyle
1;pointstyle "point";linethickness 1;lineAttributes "Solid";var1range
"1,2";num-x-gridlines 49;curveColor "[flat::RGB:0000000000]";curveStyle
"Line";rangeset"X";function \TEXUX{$(-(x-1)^{1.5}(2-x)^{1.5}+1)$};linecolor
"black";linestyle 1;pointstyle "point";linethickness 1;lineAttributes
"Solid";var1range "1,2";num-x-gridlines 49;curveColor
"[flat::RGB:0000000000]";curveStyle "Line";rangeset"X";function
\TEXUX{$((x-1)^{1.5}(2-x)^{1.5}+2)$};linecolor "black";linestyle
1;pointstyle "point";linethickness 1;lineAttributes "Solid";var1range
"1,2";num-x-gridlines 49;curveColor "[flat::RGB:0000000000]";curveStyle
"Line";rangeset"X";function \TEXUX{$(-(x-1)^{1.5}(2-x)^{1.5}+2)$};linecolor
"black";linestyle 1;pointstyle "point";linethickness 1;lineAttributes
"Solid";var1range "1,2";num-x-gridlines 49;curveColor
"[flat::RGB:0000000000]";curveStyle "Line";rangeset"X";function
\TEXUX{$\left[ 1,2,1,1\right] $};linecolor "black";linestyle 2;pointstyle
"point";linethickness 1;lineAttributes "Dash";var1range
"-5,5";num-x-gridlines 49;curveColor "[flat::RGB:0000000000]";curveStyle
"Line";function \TEXUX{$\left[ 2,2,2,1\right] $};linecolor "black";linestyle
2;pointstyle "point";linethickness 1;lineAttributes "Dash";var1range
"-5,5";num-x-gridlines 49;curveColor "[flat::RGB:0000000000]";curveStyle
"Line";function \TEXUX{$((x-1)^{1.5}(2-x)^{1.5}+2)$};linecolor
"black";linestyle 1;pointstyle "point";linethickness 1;lineAttributes
"Solid";var1range "1,2";num-x-gridlines 49;curveColor
"[flat::RGB:0000000000]";curveStyle "Line";rangeset"X";function
\TEXUX{$\left[ 1,1,2,1\right] $};linecolor "black";linestyle 1;pointplot
TRUE;pointstyle "box";linethickness 1;lineAttributes "Solid";var1range
"-5,5";num-x-gridlines 49;curveColor "[flat::RGB:0000000000]";curveStyle
"Point";function \TEXUX{$\left[ 1,2,2,2\right] $};linecolor
"black";linestyle 1;pointplot TRUE;pointstyle "box";linethickness
1;lineAttributes "Solid";var1range "-5,5";num-x-gridlines 49;curveColor
"[flat::RGB:0000000000]";curveStyle "Point";function
\TEXUX{$(-(x-3)^{1.5}(4-x)^{1.5}+1)$};linecolor "black";linestyle
1;pointstyle "point";linethickness 1;lineAttributes "Solid";var1range
"3,4";num-x-gridlines 49;curveColor "[flat::RGB:0000000000]";curveStyle
"Line";rangeset"X";function \TEXUX{$((x-3)^{1.5}(4-x)^{1.5}+2)$};linecolor
"black";linestyle 1;pointstyle "point";linethickness 1;lineAttributes
"Solid";var1range "3,4";num-x-gridlines 49;curveColor
"[flat::RGB:0000000000]";curveStyle "Line";rangeset"X";function
\TEXUX{$\left[ 3,1,3,2,4,1,4,2\right] $};linecolor "black";linestyle
1;pointplot TRUE;pointstyle "box";linethickness 1;lineAttributes
"Solid";var1range "-5,5";num-x-gridlines 49;curveColor
"[flat::RGB:0000000000]";curveStyle "Point";function \TEXUX{$\left[
3,1,4,2\right] $};linecolor "black";linestyle 1;pointstyle
"box";linethickness 1;lineAttributes "Solid";var1range
"-5,5";num-x-gridlines 49;curveColor "[flat::RGB:0000000000]";curveStyle
"Line";function \TEXUX{$\left[ 3,1,3,2,\right] $};linecolor
"black";linestyle 1;pointstyle "point";linethickness 1;lineAttributes
"Solid";var1range "-5,5";num-x-gridlines 49;curveColor
"[flat::RGB:0000000000]";curveStyle "Line";function \TEXUX{$\left[
4,1,4,2\right] $};linecolor "black";linestyle 1;pointstyle
"point";linethickness 1;lineAttributes "Solid";var1range
"-5,5";num-x-gridlines 49;curveColor "[flat::RGB:0000000000]";curveStyle
"Line";function \TEXUX{$\left[ 3,2,4,1\right] $};linecolor "black";linestyle
1;pointstyle "point";linethickness 1;lineAttributes "Solid";var1range
"-5,5";num-x-gridlines 49;curveColor "[flat::RGB:0000000000]";curveStyle
"Line";valid_file "T";tempfilename 'IIDRFK00.wmf';tempfile-properties "XPR";}%
}We consider the three cases separately.

a) For $l_{1}<l_{2}<l_{3},$ we consider first the graphs $\gamma \in \Gamma
_{L(2)}$ (on the left-hand side of the figure). By means of Lemma 3.2 in
Marinucci (2005), it is simple to ascertain that each corresponding term
produces a factor $(2l_{i}+1)$, where the $l_{i}$ correspond to the index
which is not in the two-loop (the dashed line). Direct counting of
permutations shows that there are $6$ such graphs for each fixed $l_{i}.$
Likewise, the Wigner's $6j$ coefficients arise in connection with the graphs 
$\gamma \in \Gamma _{L(3)}$ where each node is linked to all three others,
compare (\ref{6j1}) and Lemma 3.3 in Marinucci (2005); direct counting of
possible permutations shows that there can be six distinct combinations of
this form (fix for instance node 1, which by assumption is linked to all the
other three nodes: there are three degrees of freedom to choose the
connection with node 2, then two left for node 3, for a total of six, as
claimed).

b) For $l_{1}=l_{2}\neq l_{3}$, there are two possible types of graphs $%
\gamma \in \Gamma _{L(2)}$, that is, those where both 2-loops involve $l_{2}$
and those where one 2-loop involves $l_{2}$ and the other $l_{3}.$ In the
former case, there are three ways to choose the pairs, two ways in each pair
to choose the links, and two ways to choose the way to match the edges
corresponding to $l_{3}$: the total is 24. In the latter case, it can be
checked that there are 6 possible choices of pairs and 4 possible matchings
for the 2-loop which involves $l_{2}$ and $l_{3}$. The remaining term is
similar.

c) For $l_{1}=l_{2}=l_{3},$ again we choose in three possible ways the
pairs, plus the single link in two possible ways; within each pair, we can
choose $3\times 3$ couples and two possible ways to link them. The remaining
term is similar.

\hfill$\square $

\ \newline
\textbf{Remark 4.1 }Theorem 4.1 can be used to establish (in a merely \emph{%
formal} sense) Edgeworth or Cornish-Fisher type approximations for the
asymptotic behaviour of the angular bispectrum (see Hall (1991)). We provide
here only some heuristic discussion, and leave for future research the
possibility to establish rigorously a valid Edgeworth expansion. Let%
\begin{equation*}
\mu _{l_{1}l_{2}l_{3}}^{(2p)}:=E\left( \frac{I_{l_{1}l_{2}l_{3}}^{2p}}{%
\Delta _{l_{1}l_{2}l_{3}}^{p}}\right) \text{ ,}
\end{equation*}%
whence%
\begin{eqnarray*}
\sum_{p=1}^{\infty }\mu _{l_{1}l_{2}l_{3}}^{(2p)}\frac{(it)^{2p}}{2p!}
&\simeq &\sum_{p=0}^{\infty }\frac{(it)^{2p}}{2p!}(2p-1)!!+\frac{\kappa
_{4}(l_{1},l_{2},l_{3})}{24}\sum_{p=2}^{\infty }\frac{(it)^{2p}}{2p!}%
2p(2p-2)(2p-1)!!+O(l_{1}^{-2})\text{ } \\
&=&\sum_{p=0}^{\infty }\frac{(it)^{2p}}{2p!}(2p-1)!!+\frac{\kappa
_{4}(l_{1},l_{2},l_{3})}{24}\sum_{p=2}^{\infty }\frac{(it)^{2p}}{(2p-4)!}%
(2p-5)!!+O(l_{1}^{-2}) \\
&=&\sum_{p=0}^{\infty }\frac{(it)^{2p}}{2p!}(2p-1)!!+\frac{\kappa
_{4}(l_{1},l_{2},l_{3})}{24}(it)^{4}\sum_{p=0}^{\infty }\frac{(it)^{2p}}{2p!}%
(2p-1)!!+O(l_{1}^{-2}) \\
&=&\exp (\frac{-t^{2}}{2})\left\{ 1+\frac{\kappa _{4}(l_{1},l_{2},l_{3})}{24}%
t^{4}\right\} +O(l_{1}^{-2})\text{ .}
\end{eqnarray*}%
Again in a formal sense, we can take Fourier transforms on both sides,
leading to the conjecture that%
\begin{equation*}
\Pr \left\{ \frac{I_{l_{1}l_{2}l_{3}}}{\sqrt{\Delta _{l_{1}l_{2}l_{3}}}}\leq
x\right\} =\int_{-\infty }^{x}\phi (z)dz+\frac{\kappa _{4}(l_{1},l_{2},l_{3})%
}{24}\int_{-\infty }^{x}(z^{4}-6z+3)\phi (z)dz+O(l_{1}^{-2})\text{ ,}
\end{equation*}%
where $\phi (z)$ denotes the density function of a standard Gaussian
variable. It is remarkable that the first term in this expansion involves
the fourth Hermite polynomial $H_{4}(z)=(z^{4}-6z+3)$ rather than the second 
$H_{2}(z)=(z^{2}-1)$ as it is more commonly the case.

\ \newline
\textbf{Remark 4.2} By means of (\ref{formest}), it is immediate to extend
Theorem 4.1 to cover the case where the angular power spectrum is estimated
from the data; we provide below a detailed calculation for $2p=4.$ Also,
higher order moments can also be evaluated iteratively, according to the
following expression which holds for $p\geq 3$%
\begin{equation*}
EI_{l_{1}l_{2}l_{3}}^{2p}=(2p-1)!!\Delta _{l_{1}l_{2}l_{3}}^{p}\text{ }
\end{equation*}%
\begin{equation*}
+\left( 
\begin{tabular}{l}
$2p$ \\ 
$4$%
\end{tabular}%
\right) (2p-5)!!\Delta _{l_{1}l_{2}l_{3}}^{p-2}\times \left[
EI_{l_{1}l_{2}l_{3}}^{4}-3\Delta _{l_{1}l_{2}l_{3}}^{2}\right]
\end{equation*}%
\begin{equation*}
+\left( 
\begin{tabular}{l}
$2p$ \\ 
$6$%
\end{tabular}%
\right) (2p-7)!!\Delta _{l_{1}l_{2}l_{3}}^{p-3}\times \left[
EI_{l_{1}l_{2}l_{3}}^{6}-15\Delta _{l_{1}l_{2}l_{3}}^{3}\right] +...+
\end{equation*}%
\begin{equation*}
+\left( 
\begin{tabular}{l}
$2p$ \\ 
$2k$%
\end{tabular}%
\right) (2p-2k-1)!!\Delta _{l_{1}l_{2}l_{3}}^{p-k}\times \left[
EI_{l_{1}l_{2}l_{3}}^{2k}-(2k-1)!!\Delta _{l_{1}l_{2}l_{3}}^{k}\right] +...+
\end{equation*}%
\begin{equation*}
+\left( 
\begin{tabular}{l}
$2p$ \\ 
$2p-2$%
\end{tabular}%
\right) \Delta _{l_{1}l_{2}l_{3}}\left[ EI_{l_{1}l_{2}l_{3}}^{2p-2}-(2p-3)!!%
\Delta _{l_{1}l_{2}l_{3}}^{p-1}\right] +\kappa _{2p}(l_{1},l_{2},l_{3})\text{
}
\end{equation*}%
\begin{equation*}
=(2p-1)!!\Delta _{l_{1}l_{2}l_{3}}^{p}\text{ }+\sum_{k=2}^{p-1}\left( 
\begin{tabular}{l}
$2p$ \\ 
$2k$%
\end{tabular}%
\right) (2p-2k-1)!!\Delta _{l_{1}l_{2}l_{3}}^{p-k}\times \left[
EI_{l_{1}l_{2}l_{3}}^{2k}-(2k-1)!!\Delta _{l_{1}l_{2}l_{3}}^{k}\right]
+\kappa _{2p}(l_{1},l_{2},l_{3})\text{ . }
\end{equation*}%
This expression can also be exploited to derive approximations of moments,
where the term $\kappa _{2p}(l_{1},l_{2},l_{3})$ is simply neglected.

\ \ \newline
\textbf{Remark 4.3 }It is interesting to note that for $p=2$ (\ref{teo4.1})
holds with $K_{p}\equiv 0$. We provide hence an explicit evaluation of these
moments. First recall that (VMK, Chapter 9) 
\begin{eqnarray*}
\left\{ 
\begin{tabular}{lll}
$l_{1}$ & $l_{2}$ & $l_{3}$ \\ 
$l_{1}$ & $l_{2}$ & $l_{3}$%
\end{tabular}%
\right\} &=&\left[ \frac{\left\{
(l_{1}+l_{2}-l_{3})!(l_{1}+l_{3}-l_{2})!(l_{2}+l_{3}-l_{1})!\right\} ^{3}}{%
(l_{1}+l_{2}+l_{3}+1)!}\right] ^{2}\times \\
&&\sum_{n=l_{1}+l_{2}+l_{3}}^{2l_{2}+2l_{3}}\frac{(-1)^{n}(n+1)!}{%
((n-l_{1}-l_{2}-l_{3})!)^{4}((2l_{1}+2l_{2}-n)!)((2l_{1}+2l_{3}-n)!)((2l_{2}+2l_{3}-n)!)%
}
\end{eqnarray*}%
This expression can be simplified for some values of the triple $%
(l_{1},l_{2},l_{3}).$ More precisely, we have%
\begin{equation*}
\left\{ 
\begin{tabular}{lll}
$l$ & $l$ & $l$ \\ 
$l$ & $l$ & $l$%
\end{tabular}%
\right\} =\sum_{n=3l}^{4l}\frac{(-1)^{n}(n+1)!}{((n-3l)!)^{4}((4l-n)!)^{3}}
\end{equation*}%
and (VMK, eq. 9.5.2.10)%
\begin{equation*}
\left\{ 
\begin{tabular}{lll}
$l_{1}$ & $l_{2}$ & $l_{1}+l_{2}$ \\ 
$l_{1}$ & $l_{2}$ & $l_{1}+l_{2}$%
\end{tabular}%
\right\} =\frac{(2l_{1})!(2l_{2})!}{(2l_{1}+2l_{2}+1)!}\text{ .}
\end{equation*}%
Thus we have%
\begin{equation}
EI_{l_{1}l_{2},l_{1}+l_{2}}^{4}=3+\frac{6}{2l_{1}+1}+\frac{6}{2l_{2}+1}+%
\frac{6}{2l_{1}+2l_{2}+1}+6\frac{(2l_{1})!(2l_{2})!}{(2l_{1}+2l_{2}+1)!}%
\text{ , }l_{1}\neq l_{2}\text{ , }  \label{vali1}
\end{equation}%
\begin{equation}
EI_{ll,2l}^{4}=3\times 2^{2}+\frac{96}{2l+1}+\frac{24}{4l+1}+48\frac{%
((2l)!)^{2}}{(4l+1)!}\text{ ,}  \label{vali2}
\end{equation}%
and%
\begin{equation}
EI_{lll}^{4}=3\times 6^{2}+\frac{6\times 18^{2}}{2l+1}+6^{4}\frac{(l!)^{6}}{%
((3l+1)!)^{2}}\sum_{n=3l}^{4l}\frac{(-1)^{n}(n+1)!}{%
((n-3l)!)^{4}((4l-n)!)^{3}}\text{ .}  \label{vali3}
\end{equation}%
The validity of (\ref{vali1})-(\ref{vali3}) has been confirmed with a
remarkable accuracy by a Monte Carlo experiment (not reported here), where
200 replications of the sample bispectrum at various multipoles $%
(l_{1},l_{2},l_{3})$ were generated for Gaussian fields; their sample
moments were then evaluated and found to be in excellent agreement with our
theoretical results.

By means of (\ref{formest}) these results can be immediately extended to the
case where the normalization is random. More precisely, we obtain for $%
l_{1}\neq l_{2}$ 
\begin{eqnarray*}
E\widehat{I}_{l_{1}l_{2},l_{1}+l_{2}}^{4} &=&\left\{ 3+\frac{6}{2l_{1}+1}+%
\frac{6}{2l_{2}+1}+\frac{6}{2l_{1}+2l_{2}+1}+6\frac{(2l_{1})!(2l_{2})!}{%
(2l_{1}+2l_{2}+1)!}\right\} \\
&&\times \left\{ \frac{2l_{1}+1}{2l_{1}+3}\frac{2l_{2}+1}{2l_{2}+3}\frac{%
2l_{3}+1}{2l_{3}+3}\right\} \text{ ,}
\end{eqnarray*}%
\begin{eqnarray*}
E\widehat{I}_{ll,2l}^{4} &=&\left\{ 12+\frac{96}{2l+1}+\frac{24}{4l+1}+48%
\frac{((2l)!)^{2}}{(4l+1)!}\right\} \\
&&\times \left\{ \frac{2l+1}{2l+3}\frac{2l+1}{2l+5}\frac{2l+1}{2l+7}\frac{%
4l+1}{4l+3}\right\} \text{ ,}
\end{eqnarray*}%
and%
\begin{eqnarray*}
E\widehat{I}_{lll}^{4} &=&\left[ 3\times 36+\frac{6\times 18^{2}}{2l+1}+6^{4}%
\frac{(l!)^{6}}{((3l+1)!)^{2}}\sum_{n=3l}^{4l}\frac{(-1)^{n}(n+1)!}{%
((n-3l)!)^{4}((4l-n)!)^{3}}\right] \\
&&\times \prod_{k=1}^{6}\left\{ \frac{2l+1}{2l+2k-1}\right\} \text{ .}
\end{eqnarray*}

These expressions can be of practical interest for statistical inference on
cosmological data. For instance, the square bispectrum is often used in
goodness-of-fit statistics to test the validity of the Gaussian assumption.
The previous results yield immediately its exact variance, which so far has
been typically evaluated by Monte Carlo simulations (see for instance
Komatsu et al. (2002)).

\ 

\begin{center}
\textbf{REFERENCES}
\end{center}

\qquad

\textbf{Arjunwadkar, M., C.R. Genovese, C.J. Miller, R.C. Nichol and L.
Wasserman (2004) }``Nonparametric Inference for the Cosmic Microwave
Background'', \emph{Statistical Science, }Vol. 19, Issue 2, pp.308-321

\textbf{Bartolo, N., S. Matarrese and A. Riotto (2002) }``Non-Gaussianity
from Inflation'', \emph{Physical Review D}, Vol.65, Issue 10, id. 3505; also
available at http://it.arxiv.org as astro-ph/0112261

\textbf{Babich, D. (2005) }``Optimal Estimation of Non-Gaussianity'',
preprint, available at http://it.arxiv.org as astro-ph/0503375

\textbf{Biedenharn, L.C. and J.D. Louck (1981) }\emph{The Racah-Wigner
Algebra in Quantum Theory, }Encyclopedia of Mathematics and its
Applications, Volume 9, Addison-Wesley

\textbf{Dor\`{e}, O., S. Colombi, F.R. Bouchet (2003)} ``Probing CMB
Non-Gaussianity Using Local Curvature'', \emph{Monthly Notices of the Royal
Astronomical Society}, Vol. 344, Issue 3, pp. 905-916, available at
http://it.arxiv.org as astro-ph/0202135

\textbf{Foulds, L.R. (1992) }\emph{Graph Theory Applications, }%
Springer-Verlag

\textbf{Giraitis, L. and D. Surgailis (1987) }``Multivariate Appell
Polynomials and the Central Limit Theorem'', in \emph{Dependence in
Probability and Statistics}, Birkhauser, pp. 21--71

\textbf{Hall, P. (1991) }\emph{The Bootstrap and Edgeworth Expansion},
Springer-Verlag

\textbf{Hansen, F.K., D. Marinucci and N. Vittorio (2003) }``The Extended
Empirical Process Test for Non-Gaussianity in the CMB, with an Application
to Non-Gaussian Inflationary Models'', \emph{Physical Review D,} Vol. 67,
Issue 12, id. 3004; also available at http://it.arxiv.org as astro-ph/0302202

\textbf{Hu, W. (2001) }``The Angular Trispectrum of the CMB'', \emph{%
Physical Review D}, Vol. 64, Issue 8, id.3005; also available at
http://it.arxiv.org as astro-ph/0105117

\textbf{Jin, J., J.L. Starck, D. Donoho, \ N. Aghanim and O. Forni (2004)}
``Cosmological non-Gaussian Signature Detection: Comparing the Performance
of Different Statistical Tests'', \emph{Eurasip Journal on Applied Signal
Processing}, forthcoming

\textbf{Kim, P.T. and J.-Y. Koo (2002) }``Optimal Spherical Deconvolution'', 
\emph{Journal of Multivariate Analysis}, Vol. 80, Issue 1, pp 21-42

\textbf{Kim, P.T., Koo, J.-Y. and H.J. Park (2004) }``Sharp Minimaxity and
Spherical Deconvolution for Super-Smooth Error Distributions'', \emph{%
Journal of Multivariate Analysis}, Vol. 90, Issue 2, pp. 384-392

\textbf{Komatsu, E. and D.N. Spergel (2001) }``Acoustic Signatures in the
Primary Microwave Background Bispectrum'', \emph{Physical Review D},\emph{\ }%
Vol. 63, Issue 6, id. 3002; also available at http://it.arxiv.org as
astro-ph/0005036

\textbf{Komatsu, E. et al. (2002) }``Measurement of the Cosmic Microwave
Background Bispectrum on the COBE DMR Sky Maps'', \emph{Astrophysical
Journal, }Vol. 566, pp.19-29

\textbf{Komatsu, E. et al. (2003) }``First Year \emph{Wilkinson Microwave
Anisotropy Probe (WMAP)} Observations: Tests of Gaussianity'', \emph{%
Astrophysical Journal Supplement Series}, Vol. 148, Issue 1, pp.119-134;
also available at http://it.arxiv.org as astro-ph/0302223

\textbf{Leonenko, N.N. (1999) }\emph{Limit Theorems for Random Fields with
Singular Spectrum, }Kluwer, Dordrecht

\textbf{Loh, W.L. (2005)} ``Fixed-domain Asymptotics for a Subclass of
Matern-type Gaussian Random Fields'', \emph{Annals of Statistics}, to appear.

\textbf{Marinucci, D. (2004) }``Testing for non-Gaussianity on Cosmic
Microwave Background Radiation: a Review'', \emph{Statistical Science, }Vol.
19, Issue 2, pp.294-307

\textbf{Marinucci, D. and M. Piccioni (2004) }``The Empirical Process on
Gaussian Spherical Harmonics'', \emph{Annals of Statistics, }Vol. 32, Issue
3, pp.1261-1288.

\textbf{Marinucci, D. (2005) }``High Resolution Asymptotics for the Angular
Bispectrum'', \emph{Annals of Statistics, }forthcoming

\textbf{Park, C.-G. (2004) }``Non-Gaussian Signatures in the Temperature
Fluctuation Observed by the Wilkinson Microwave Anisotropy Probe'', \emph{%
Monthly Notices of the Royal Astronomical Society, }Vol.349, Issue 1,
p.313-320

\textbf{Peacock, J.A. (1999) }\emph{Cosmological Physics, }Cambridge
University Press, Cambridge.

\textbf{Peebles, P.J.E. (1993) }\emph{Principles of Physical Cosmology, }%
Princeton University Press, Princeton.

\textbf{Phillips, N.G. and A. Kogut (2000) }``Statistical Power, the
Bispectrum and the Search for Non-Gaussianity in the CMB Anisotropy'', \emph{%
Astrophysical Journal, }Vol. 548, Issue 2, pp. 540-549;\emph{\ }also
available at http://it.arxiv.org as astro-ph/0010333

\textbf{Stein, M.L. (1999) }\emph{Interpolation of Spatial Data. Some Theory
for Kriging. }Springer-Verlag.

\textbf{Varshalovich, D.A., A.N. Moskalev, and V.K. Khersonskii (1988), }%
\emph{Quantum Theory of Angular Momentum, }World Scientific, Singapore

\textbf{Vilenkin, N.J. and A.U. Klimyk (1991) }\emph{Representation of Lie
Groups and Special Functions, }Kluwer, Dordrecht

\textbf{Yaglom, A.M. (1986) }\emph{Correlated Theory of Stationary and
Related Random Functions I. Basic Results, }Springer-Verlag.\newline

\ 

Address for correspondence:

Dipartimento di Matematica, Universita' di Roma Tor Vergata

via della Ricerca Scientifica 1, Roma, Italy. Postal Code: 00133

\end{document}